\documentclass[a4paper,oneside,11pt,final]{article}

\usepackage{geometry}
\usepackage{array}
\usepackage{amssymb,amsmath,amsfonts,amsthm}
\usepackage{showkeys}
\usepackage{color}

\usepackage{tikz,pgfplots}

\if@onethmnum
  \newtheorem{theorem}{Theorem}
  \newtheorem{lemma}[theorem]{Lemma}
  \newtheorem{corollary}[theorem]{Corollary}
  \newtheorem{proposition}[theorem]{Proposition}
  \newtheorem{definition}[theorem]{Definition}
\else
  \newtheorem{theorem}{Theorem}
  \newtheorem{lemma}[theorem]{Lemma}
  \newtheorem{corollary}[theorem]{Corollary}
  \newtheorem{proposition}[theorem]{Proposition}
  \newtheorem{definition}[theorem]{Definition}
\fi

\if@onethmnum
  \newtheorem{conjecture}{Conjecture}
  \newtheorem{example}{Example}
 \else

\fi

\if@onethmnum
  \newtheorem{assumption}{Assumption}
\else
  \newtheorem{assumption}{Assumption}
\fi

\title{Projection methods and discrete gradient methods for preserving first integrals of ODEs}
\author{R. A. Norton, D. I. McLaren, G. R. W. Quispel, A. Stern \& A. Zanna}

\newcommand{\bbR}{\mathbb{R}}
\newcommand{\wS}{\widetilde{S}}

\newcommand{\tff}{\tilde{f}}
\newcommand{\tgg}{\tilde{g}}
\newcommand{\tii}{\tilde{i}}
\newcommand{\tjj}{\tilde{j}}
\newcommand{\bii}{\bar{i}}
\newcommand{\bjj}{\bar{j}}
\newcommand{\hii}{\hat{i}}
\newcommand{\brii}{\breve{i}}
\newcommand{\tsfrac}  [2] { {\textstyle \frac{#1}{#2} } }

\newcommand{\fourbyone}[4]{\begin{array}{c} #1 \\ #2 \\ #3 \\ #4 \end{array} }
\newcommand{\dd}{{\rm d}}

\newcommand{\opspan}{\operatorname{span}}
\newcommand{\oprange}{\operatorname{range}}
\newcommand{\opnull}{\operatorname{null}}

\newcommand{\cA}{\mathcal{A}}
\newcommand{\cB}{\mathcal{B}}
\newcommand{\cS}{\mathcal{S}}
\newcommand{\cM}{\mathcal{M}}

\newcommand*{\matminus}{%
  \leavevmode
  \hphantom{0}%
  \llap{%
    \settowidth{\dimen0 }{$0$}%
    \resizebox{1.1\dimen0 }{\height}{$-$}%
  }%
}

\begin{document}
\maketitle

\begin{abstract}
In this paper we study linear projection methods for approximating the solution and simultaneously preserving first integrals of autonomous ordinary differential equations.  We show that (linear) projection methods are a subset of discrete gradient methods.  In particular, each projection method is equivalent to a class of discrete gradient methods (where the choice of discrete gradient is arbitrary) and earlier results for discrete gradient methods also apply to projection methods.  Thus we prove that for the case of preserving one first integral, under certain mild conditions, the numerical solution for a projection method exists and is locally unique, and preserves the order of accuracy of the underlying method.  

In the case of preserving multiple first integrals the relationship between projection methods and discrete gradient methods persists.  Moreover, numerical examples show that similar existence and order results should also hold for the multiple integral case.

For completeness we show how existing projection methods from the literature fit into our general framework.
\end{abstract}


\section{Introduction}
\label{sec intro}

First, consider an autonomous ordinary differential equation (ODE) with only one first integral.  We will consider the case of multiple first integrals later.  We consider the same problem as in \cite{nortonquispel1}:
\begin{equation}
\label{p1}
	\dot{x} = f(x) \qquad t > 0, 
\end{equation}
where $x(t) \in \bbR^d$ for some $d \in \mathbb{N}$, $x(0) = x_0 \in \bbR^d$ is the initial condition and $f:\bbR^d \rightarrow \bbR^d$ is locally Lipschitz continuous.  Existence theory for ODEs (see e.g. \cite[Thm. I.7.3 on p.37]{SolvingODEs1}) implies that given a bounded set $B \subset \bbR^d$, there exists a $T > 0$ such that for any $x_0 \in B$ the solution exists and remains bounded for $t \in [0,T]$.  We assume that \eqref{p1} has a conserved first integral $I : \bbR^d \rightarrow \bbR$ so that 
$$
	x(t) \in \cM	_{x_0}:= \{ z \in \bbR^d : I(z) = I(x_0) \} \qquad \mbox{for all $t \in [0,T]$.}
$$

To simplify the notation define $i := \nabla I$ and let us also assume that $I$ is a Morse function (i.e. smooth with non-degenerate critical points) and that $i : \bbR^d \rightarrow \bbR^d$ is locally Lipschitz continuous.  As in \cite{MQR99}, and discussed in detail in \cite{nortonquispel1}, if $i(x(t)) \neq 0$ for $t > 0$ then we may write \eqref{p1} as
\begin{equation}
\label{p3}
	\dot{x} = S(x) i(x)
\end{equation}
where $S: \bbR^d \rightarrow \bbR^{d \times d}$ is a skew-symmetric ($S^T = -S$) matrix-valued function.   In general, $S$ is not unique.  One choice for $S$ is the so-called \emph{default} formula,
\begin{equation}
\label{eqn3}
	S(x) = \frac{f(x) i(x)^T - i(x) f(x)^T}{|i(x)|^2}.
\end{equation}
Since $I$ is a Morse function, the default $S$ is locally bounded on $\{ x \in \bbR^d : i(x) \neq 0 \}$ and for a bounded set $B \subset \bbR^d$ there exists a constant $C_1 = C_1(B)$ such that 
\begin{equation}
\label{eqn2}
	|f(x)| \leq C_1 |i(x)| \qquad \mbox{for all $x \in B$}.
\end{equation}
Also define $C_2 = C_2(B) := C_1 + \tsfrac{1}{5}$.  

In general it will be beneficial to approximate the solution to \eqref{p1} in such a way so that $I$ is preserved exactly (in practice up to round off error or a specified tolerance) by the approximate solution.  Both (linear) projection methods (see e.g. \cite[\S IV.4 and \S V.4.1]{HLW} and references therein) and discrete gradient methods (see e.g. \cite{MQR99,quispelcapel96,quispelturner96}) are types of methods that achieve this.  In the special case of Hamiltonian systems one must choose whether to preserve the Hamiltonian integral or the symplectic structure (only the exact solution up to time rescaling preserves both, see e.g. \cite{ge88}), but there are many examples where preserving the Hamiltonian is advantageous, see e.g. \cite{quispelmclaren08,simo92}.

First, let us define linear projection methods.  The basic idea of a projection method is to couple a one-step method with a projection so that after a full time step the approximate solution to the ODE lies on the manifold $\cM_{x_0}$.  Let $\tff : \bbR^d \times \bbR^d \times [0,\infty) \rightarrow \bbR^d$ define an arbitrary one-step method applied to \eqref{p1} with time step $h$, so that 
\begin{equation}
\label{e1}
	\frac{x'-x}{h} = \tff(x,x',h)
\end{equation}
where $x = x_n$ and $x' = x_{n+1}$ at each time step\footnote{It will be our convention to let $x = x_n$ (the approximate solution at step $n$) and $x' = x_{n+1}$.}.  Then, one step of a linear projection method (c.f. \cite[Algorithm IV.4.2]{HLW}) $x \mapsto x'$ is defined by:  Given $x \in \bbR^d$ and $h \in [0,\infty)$,
\begin{enumerate}
\item compute $y$ such that $y = x + h \tff(x,y,h)$,
\item compute $x' \in \mathcal{M}_x$ by projecting $y$ onto $\mathcal{M}_x$.
\end{enumerate}
In this paper we will only concern ourselves with linear projections so that step 2 of the above algorithm is given by: 
\begin{enumerate}
\setcounter{enumi}{1}
\item
compute $x' \in \cM_x$ by solving $x' = y + \lambda \tii(x,x',h)$ and $I(x') = I(x)$ for $x' \in \bbR^d$ and $\lambda \in \bbR$,
\end{enumerate}
where $\tii : \bbR^d \times \bbR^d \times [0,\infty) \rightarrow \bbR^d$ is a vector field that defines the direction of the projection and is typically an approximation of $i$.  We refer to this type of projection as a \emph{linear} projection because $(x' - y) \parallel \tii(x,x',h)$.

Note that for a method defined by $y = x + h \tff(x,y,h)$ in step 1 of the algorithm above, there exists an implicitly defined map $\Phi_h : \bbR^d \rightarrow \bbR^d$ such that $y = \Phi_h(x)$.  If we define $\tgg(x,h) := (\Phi_h(x)-x)/h$, then we may alternatively write step 1 as $y = x + h \tgg(x,h)$.  Using $\tgg$ instead of $\tff$ in step 1 allows us to easily eliminate $y$ from the algorithm and express the algorithm in a single line: Given $x \in \bbR^d$ and $h \in [0,\infty)$ compute $x' \in \bbR^d$ and $\lambda \in \bbR$ such that
$$
	x' = x + h \tgg(x,h) + \lambda \tii(x,x',h) \qquad \mbox{and} \qquad I(x') = I(x).
$$
For more generality in our projection methods, in addition to allowing different choices of $\tii$, we will also modify $\tgg$ so that it may also depend on $x'$.  Switching back to using $\tff$ instead of $\tgg$ in the notation we get our general form of a linear projection method for preserving a single first integral: Given $x \in \bbR^d$ and $h \in [0,\infty)$ compute $x' \in \bbR^d$ and $\lambda \in \bbR$ such that 
\begin{equation}
\label{eqn0}
	x' = x + h \tff(x,x',h) + \lambda \tii(x,x',h) \qquad \mbox{and} \qquad  I(x') = I(x).
\end{equation}
By choosing $\tff$ and $\tii$ differently, we obtain different projection methods.  To the best of our knowledge all of the projection methods that have been described in the literature fit into this framework (we are only aware of linear projection methods but it may be possible to define projection methods in spaces that are not linear spaces), including the (non-symmetric) standard projection method in \cite[Algorithm IV.4.2]{HLW} and the symmetric projection method in \cite[\S V.4.1]{HLW} and \cite{hairer00}.   This will be discussed in more detail in Section \ref{sec existing}.

The other type of integral preserving methods we consider are discrete gradient methods.  For their definition we must first define a \emph{discrete gradient of $I$} - a special type of discretization of the gradient of $I$.

\begin{definition}
\emph{(Gonzalez \cite{gonzalez96})}
A discrete gradient of $I$, denoted $\bii: \bbR^d \times \bbR^d \rightarrow \bbR^d$, is continuous and satisfies
$$
	\bii(x,x') \cdot (x' - x) = I(x') - I(x) \quad \mbox{and} \quad \bii(x,x) = i(x) \quad \mbox{for all $x,x' \in \bbR^d$.}
$$
\end{definition}

Formulae for constructing discrete gradients include the one used in the average (or averaged) vector field method (called \emph{mean value discrete gradient} in \cite{MQR99}, also see \cite{quispelmclaren08}) and the coordinate increment method \cite{abeitoh}.

If we let $\bii$ be a discrete gradient of $I$ and $\wS : \bbR^d \times \bbR^d \times [0,\infty) \rightarrow \bbR^{d \times d}$ be a skew symmetric continuous and consistent approximation of $S$ then a discrete gradient method for solving \eqref{p1} is defined by the mapping $x \mapsto x'$ where
\begin{equation}
\label{disc1}
	x' = \begin{cases}
		x + h \wS(x,x',h) \bii(x,x') & \mbox{if $i(x) \neq 0$}, \\
		x & \mbox{if $i(x) = 0$}.
	\end{cases}
\end{equation}
In this paper we only consider the large class of discrete gradients where $\wS$ is defined by the formula 
\begin{equation}
\label{eqn4}
	\wS(x,x',h) := \frac{ \tff(x,x',h) \tii(x,x',h)^T - \tii(x,x',h) \tff(x,x',h)^T}{\hii(x,x',h) \cdot \brii(x,x',h)},
\end{equation}
where $\tff : \bbR^d \times \bbR^d \times [0,\infty) \rightarrow \bbR^d$ is a continuous consistent approximation of $f$ and $\tii$, $\hii$ and $\brii$ are all maps from $\bbR^d \times \bbR^d \times [0,\infty)$ to $\bbR^d$ and are continuous consistent approximations of $i$.  

All discrete gradient methods of this type preserve $I$ because
$$
	I(x') - I(x) = \bii(x,x') \cdot (x'-x) = h (\bii(x,x'))^T \wS(x,x',h) \bii(x,x') = 0
$$
for all $h$, $x$ and $x'$ satisfying \eqref{disc1}, where $\bii$ is a discrete gradient of $I$.  The final equality follows from the fact that $\wS$ is skew symmetric.

In \cite{nortonquispel1} discrete gradient methods of this type were studied and it was shown that under certain local Lipschitz conditions and for sufficiently small time step the numerical solution to \eqref{disc1} (with $\wS$ defined by \eqref{eqn4}) exists and is locally unique, independent of the distance to critical points of $I$.  For arbitrary $p \in \mathbb{N}$ it was also shown how to construct discrete gradient methods that have order of accuracy $p$.

In this paper we will show that all linear projection methods of the type \eqref{eqn0} are equivalent to discrete gradient methods where the approximate solution is independent of the particular choice of discrete gradient $\bii$.  ..We prove this by showing that each projection method is equivalent to (generally) several discrete gradient methods, in the sense that a projection method and several discrete gradient methods (defined with different choices for $\bii$) all define the same map $x \mapsto x'$ for a given $h$.  A consequence of this result is that projection methods are a subset of discrete gradient methods. 

In this case when there is only one first integral to preserve, we can then use the theory in \cite{nortonquispel1} to obtain by simple corollary new existence, uniqueness and order of accuracy results for a large number of linear projection methods (only restricted by certain mild local Lipschitz conditions on $\tff$ and $\tii$).

When there is more than one first integral to preserve, we will prove that the same equivalence between discrete gradient and linear projection methods holds, and as a consequence projection methods are a subset of discrete gradient methods for the multiple integral situation.  Since the theory in \cite{nortonquispel1} is only for the single first integral case we do not obtain new results about existence, uniqueness and order of accuracy from discrete gradient method theory for the multiple first integral case.  Proving these results for general linear projection methods and discrete gradient methods in the multiple integral case remains an open problem.

The remainder of this paper is organised as follows.  In Section \ref{sec equiv1} we prove our first result about the equivalence of linear projection methods and a class of discrete gradient methods in the case when \eqref{p1} has a single first integral.  Then, in Section \ref{sec results} we use this result and the theory from \cite{nortonquispel1} to get new results about existence, local uniqueness, and order of accuracy for linear projection methods.  In Section \ref{sec existing} we demonstrate how several projection methods already described in the literature are special cases in our framework and how our new results improve on existing results by allowing more freedom on the projection direction than previously, and our results are independent of the distance to critical points of $I$.  In Section \ref{sec equiv2} we then consider the case when \eqref{p1} has more than one first integral and our projection and discrete gradient methods are designed to preserve multiple first integrals.  We write down a new expression for linear projection methods in this case involving oblique projection matrices and prove equivalence with discrete gradient methods.  Using numerical experiments we illustrate how the order of accuracy results, proven in the single first integral case, also appear to hold in the multiple integral case.  Finally, in Section \ref{sec conclusion} we discuss the implications of this work and possible avenues for future research.  

\section{Equivalence in the single first integral case}
\label{sec equiv1}

In this section we explore the relationship between linear projection methods and discrete gradient methods.  We will see that each linear projection method is equivalent to possibly several discrete gradient methods where the choice of discrete gradient is arbitrary.  Note, however, that discrete gradient methods are not always projection methods so that projection methods are a subset of discrete gradient methods.

So far we have not proven that the projection method defined by \eqref{eqn0} is well-defined in the sense that the implicit system of equations \eqref{eqn0} for $x'$ and $\lambda$ has a unique solution for sufficiently small time step $h$.  So let us assume that $h$ is sufficiently small and $x'$ and $\lambda$ are uniquely defined by \eqref{eqn0} and that $\tii \cdot \bii \neq 0$ (definitions of $\tii$ and $\bii$ will be given).  In the next section we will prove an existence result that justifies these two assumptions under sufficient conditions for $\tff$, $\tii$, $\bii$ and $h$. 

The following theorem shows that linear projection methods may be expressed in several equivalent ways and the following corollary explains how each linear projection method is equivalent to possibly several discrete gradient methods where the choice of discrete gradient method is arbitrary.

\begin{theorem}
\label{thm equiv1}
Let $\tii : \bbR^d \times \bbR^d \times [0,\infty) \rightarrow \bbR^d$ be a consistent approximation of $i$, let $\bii$ be an arbitrary discrete gradient of $I$, and let $\tff : \bbR^d \times \bbR^d \times [0,\infty) \rightarrow \bbR^d$ be a consistent approximation of $f$.

Assuming that given $x \in \bbR^d$ and $h \in [0,\infty)$, each of the methods below have uniquely defined $x'$ and $\lambda$, and that $\tii \cdot \bii \neq 0$, then they define the same linear projection method.
\begin{align}
	x' &= x + h \tff(x,x',h) + \lambda \tii(x,x',h) &&\mbox{such that } I(x') = I(x), \label{eqn1st} \\
	x' &= x + h P(x,x',h) \tff(x,x',h) && \mbox{where  } P(x,x',h) = I - \tsfrac{\tii(x,x',h) \, \bii(x,x')^T}{\tii(x,x',h)^T \bii(x,x')}, \label{eqn2nd} \\
	x' &= x + h \wS(x,x',h) \bii(x,x') && \mbox{where } {\scriptstyle \wS(x,x',h) = \frac{\tff(x,x',h) \, \tii(x,x',h)^T - \tii(x,x',h) \tff(x,x',h)^T}{\tii(x,x',h) \cdot \bii(x,x',h)}}. \label{eqn3rd}
\end{align}
\end{theorem}

\begin{proof}
For given $x \in \bbR^d$ and $h \in [0,\infty)$ suppose $x' \in \bbR^d$ and $\lambda \in \bbR$ satisfy \eqref{eqn1st}.  Take the inner product of \eqref{eqn1st} with $\bii(x,x')$ to get
$$
	0 = I(x') - I(x) = \bii \cdot (x' - x) = h \tff \cdot \bii + \lambda \tii \cdot \bii, \qquad \mbox{hence} \qquad \lambda = - h \frac{\tff \cdot \bii}{\tii \cdot \bii}.
$$
Substituting this into \eqref{eqn1st} we find that
$$
	x' = x + h \tff - h \frac{\tff\cdot\bii}{\tii\cdot\bii} \tii = x + h \left( \tff - \frac{\tii \, \bii^T}{\tii^T \bii} \tff \right) = x + h P \tff,
$$
and
$$ 
	x' = x + h P\tff = x + h \left( \frac{\tii^T \bii}{\tii^T \bii} \tff - \frac{\tii \, \bii^T}{\tii^T \bii} \tff \right) 
	= x + h \left( \frac{\tff \, \tii^T }{\tii^T \bii} \bii - \frac{\tii  \tff^T }{\tii^T \bii} \bii \right)
	= x + h \wS \bii.
$$
Thus, if $x'$ and $\lambda$ satisfy \eqref{eqn1st} then $x'$ satisfies \eqref{eqn2nd} and \eqref{eqn3rd}.  To see the converse, note that if $x'$ satisfies \eqref{eqn2nd} or \eqref{eqn3rd} then $I(x') = I(x)$ (take the inner product of \eqref{eqn2nd} or \eqref{eqn3rd} with $\bii(x,x')$ and use the definition of a discrete gradient).  Define $\lambda := -h (\tff \cdot \bii)/(\tii\cdot \bii)$, then $x'$ and $\lambda$ satisfy \eqref{eqn1st}.
\end{proof}

Note that \eqref{eqn1st} is the same as \eqref{eqn0}, our general form for a linear projection method.  In \eqref{eqn2nd}, $P$ is a projection matrix satisfying $P \tff \perp \bii$ and $(I-P) \tff \parallel \tii$, i.e. the range of $P$ is $\operatorname{span} \lbrace \bii \rbrace^\perp$ and the null space of $P$ is $\operatorname{span}\lbrace \tii \rbrace$.

Using the equivalence of \eqref{eqn1st} and \eqref{eqn3rd} we get the following corollary.

\begin{corollary}
\label{cor equiv}
Let $\tii : \bbR^d \times \bbR^d \times [0,\infty) \rightarrow \bbR^d$ be a consistent approximation of $i$, let $\tff : \bbR^d \times \bbR^d \times [0,\infty) \rightarrow \bbR^d$ be a consistent approximation of $f$ and let $\bii$ be an arbitrary discrete gradient of $I$.  If we define $\hii \equiv \tii$ and $\brii \equiv \bii$ (or vice versa), then the linear projection method defined by \eqref{eqn0} is equivalent to the discrete gradient method defined by \eqref{disc1} where $\wS$ is defined by \eqref{eqn4}.
\end{corollary}

Defining $\hii \equiv \tii$ and $\brii \equiv \bii$ in the definition of a discrete gradient method is a restriction so linear projection methods are a subset of all possible discrete gradient methods.  

Also notice that the methods described by \eqref{eqn2nd} and \eqref{eqn3rd} depend on an arbitrarily chosen discrete gradient $\bii$, whereas linear projection methods are independent of $\bii$.  At first glance it would appear that the mapping $x \mapsto x'$ defined by \eqref{eqn2nd} and \eqref{eqn3rd} should depend on the choice of $\bii$ and these methods would give different approximations to \eqref{p1} for different choices of $\bii$.  It is perhaps surprising that this is not the case, and (as a consequence of Theorem \ref{thm equiv1} since \eqref{eqn1st} does not depend on $\bii$) they give the same approximation to \eqref{p1} regardless of how $\bii$ is chosen.  Thus, each linear projection method defines an equivalence class of discrete gradient methods, and is uniquely defined by choosing $\tff$ and the direction of projection given by $\tii$.

\section{Existence, uniqueness and order of accuracy}
\label{sec results}

In this section we will exploit the equivalence between linear projection methods and discrete gradient methods by using theory developed for discrete gradient methods to prove new results about linear projection methods.  

Typically, the projection step of a projection method (step 2 in our original algorithm) requires solving an implicit nonlinear system of equations, and a new system of equations must be solved at each time step.  A basic question regarding projection methods is: Does there exist a unique solution to each of these systems of equations?  A further question is: Does a projection method retain the same order of accuracy as the underlying method (the underlying method is step 1 in our original algorithm)?  

Linear projection methods have already been studied in the literature (see e.g. \cite[\S IV.4 and \S V.4.1]{HLW} and \cite{hairer00}) and questions of existence and uniqueness, and order of accuracy have already been answered in some cases.  However, these results were only stated for particular special cases of $\tii$ (see Section \ref{sec existing}) and their proofs rely on either a simple geometric argument (which only holds for the standard projection method when $\tii(x,x',h) := i(x')$), the Implicit Function Theorem, or the Newton-Kantorovich Theorem.  Closer examination of these techniques -- with the assistance of results in \cite{papi05} and \cite{ortega} that give a lower bound on the radius of existence for the Implicit Function Theorem and the Newton-Kantorovich Theorem -- reveals that the time step restriction on $h$ (or radius of existence) for existence of the numerical solution is $h \leq C |i(x)|^r$ for some positive constants $C$ and $r$.  If $x$ is near to a critical point of $I$ (so that $i(x) \approx 0$) then this type of restriction is undesirable and in numerical simulations it appears to be unnecessary.  Our new results below are an improvement and extension on these earlier results because we avoid this restriction, and we only place mild Lipschitz continuity conditions on the projection direction $\tii$ so that the results hold for a much wider class of projection methods.

For the following results we require the following definition of a ball around a point $x \in \bbR^d$.  Given $x \in \bbR^d$ and a constant $R > 0$ define 
$$
	B_R(x) := \left\{ z \in \bbR^d : |z-x| \leq \tsfrac{|i(x)|}{R} \right\}.
$$ 
Note that when $i(x) = 0$ we have $B_R(x) = \{ x \}$.

To simplify the presentation that follows let us define several `Assumptions'.

\begin{assumption}
\label{assumption1} 
Given a bounded set $B \subset \bbR^d$, we say that $\tff : \bbR^d \times \bbR^d \times [0,\infty) \rightarrow \bbR^d$ satisfies Assumption \ref{assumption1} for positive constants $R$, $L$ and $H$ if
\begin{align*}
	\tff(x,x,0) &= f(x) \\
	|\tff(u,v,h) - \tff(w,v,h)| &\leq L |u-w| \\
	|\tff(u,v,h) - \tff(u,w,h)| &\leq L |v - w| \\
	|\tff(x,x,h) - \tff(x,x,0)| &\leq L h |i(x)|,
\end{align*}
for all $u,v,w \in B_R(x)$, $x \in B$ and $h \in [0,H)$.
\end{assumption}

\begin{assumption}
\label{assumption2}
Given a bounded set $B \subset \bbR^d$, we say that $\tii : \bbR^d \times \bbR^d \times [0,\infty) \rightarrow \bbR^d$ satisfies Assumption \ref{assumption2} for positive constants $R$, $L$ and $H$ if
\begin{align*}
	\tii(x,x,0) &= i(x) \\
	|\tii(u,v,h) - \tii(w,v,h)| &\leq L |u-w| \\
	|\tii(u,v,h) - \tii(u,w,h)| &\leq L |v - w| \\
	|\tii(x,x,h) - \tii(x,x,0)| &\leq L h |i(x)|,
\end{align*}
for all $u,v,w \in B_R(x)$, $x \in B$ and $h \in [0,H)$.

We say that $\bii : \bbR^d \times \bbR^d \rightarrow \bbR^d$ satisfies Assumption \ref{assumption2} for $R$ and $L$ if $\tii : \bbR^d \times \bbR^d \times [0,\infty) \rightarrow \bbR^d$ defined by $\tii(x,x',h) = \bii(x,x')$ for all $x,x' \in \bbR^d$ and $h \in [0,\infty)$ satisfies Assumption \ref{assumption2} for $R$, $L$ and any positive constant $H$.
\end{assumption}

\begin{assumption}
\label{assumption3}
Given a bounded set $B \subset \bbR^d$, we say that $g : \bbR^d \rightarrow \bbR^d$ satisfies Assumption \ref{assumption3} for positive constants $R$ and $L$ if 
$$
	|g(u) - g(v)| \leq L |u-v|
$$
for all $u,v \in B_R(x)$ and all $x \in B$.
\end{assumption}

Note that since we have assumed that $f$ is locally Lipschitz continuous, for any $R > 0$ there exists a corresponding $L > 0$ such that $f$ satisfies Assumption \ref{assumption3}.  Similarly for $i$.

The following theorem ensures for sufficiently small $h$ and under certain local Lipschitz continuity conditions, that linear projection methods (defined by \eqref{eqn0}) have a numerical solution that is locally unique.  Its proof is omitted because it is a direct consequence of Theorem 2.1 in \cite{nortonquispel1} and Corollary \ref{cor equiv} where $\bii$ is chosen to be an arbitrary discrete gradient of $I$ satisfying Assumption \ref{assumption2} for $R$, $L$ and $H$ defined as in the theorem below.

\begin{theorem}
\label{thm exist}
Let $B$ be a bounded set in $\bbR^d$, let $C_2$ be the constant defined by \eqref{eqn2}, and suppose that $R$, $L$ and $H$ are positive constants such that 
\begin{enumerate}
\item $\tff : \bbR^d \times \bbR^d \times [0,\infty) \rightarrow \bbR^d$ satisfies Assumption \ref{assumption1} for $R$, $L$ and $H$, and
\item $\tii: \bbR^d \times \bbR^d \times [0,\infty) \rightarrow \bbR^d$ satisfies Assumption \ref{assumption2} for $R$, $L$ and $H$.  
\end{enumerate}
Define
$$
	R' := \max \lbrace R, 10 L \rbrace \quad \mbox{and} \quad 
	H' := \min \left\{ H, \tsfrac{1}{10L} , \tsfrac{1}{6C_2 R'},\tsfrac{1}{(36C_2 + 6) L} \right\}.
$$
Then for each $x \in B$ and $h \in [0,H')$, there exists a unique $x' \in B_{R'}(x)$ and $\lambda \in \bbR$ such that $|\lambda| \leq \tsfrac{11}{5} C_2 h$ satisfying \eqref{eqn0}.  
\end{theorem}

This existence result only provides us with local uniqueness since we are only sure that $x'$ is unique in the ball $B_{R'}(x) \subset \bbR^d$.  

Now let us consider the order of accuracy of linear projection methods.  We use the following definition for order of accuracy, which is similar to \cite[Def. V.1.3]{gautschi}.

\begin{definition}
A one-step method $x \mapsto x'$ with time step $h$ for solving \eqref{p1} has order of accuracy $p \in \mathbb{N}$, if for problems with sufficiently smooth $f$ there exist positive constants $C$ and $H$ such that
$$
	|x' - x(t+h)| \leq C h^{p+1} \qquad \mbox{for all $h \in [0,H]$ and all $x \in B$,}
$$
where $x(\cdot)$ denotes the solution to \eqref{p1} with $x(t) = x$ for some $t\geq0$ and $B$ is a compact set in $\bbR^d$.  The constants $C$ and $H$ may depend on $B$ but should be independent of $x$ and $h$.  
\end{definition}

If we are given an underlying method that is of order $p$ for some $p \in \mathbb{N}$, i.e. the method $x \mapsto y$ defined by $y = x + h \tff(x,y,h)$ is of order $p$, then an important question to ask is:  What additional conditions (in addition to Assumption \ref{assumption2}) on $\tii$ (recall $\tii$ defines the direction of the projection) are required to ensure that a linear projection method defined by \eqref{eqn0} is also of order $p$?  The following theorem gives the answer: none!  Besides Assumption \ref{assumption2}, there are no additional conditions on $\tii$ that are required to ensure a linear projection method is of order $p$.

Again, we rely on theory in \cite{nortonquispel1} to achieve our result.  The following theorem is a special case of Theorem 3.3 in \cite{nortonquispel1} using Corollary \ref{cor equiv} and an arbitrary discrete gradient $\bii$ satisfying Assumption \ref{assumption2}.

\begin{theorem}
\label{thm order}
For a compact set $B \subset \bbR^d$, let $C_2$, $R$, $L$, $H$, $\tff$, $\tii$, $R'$ and $H'$ be defined as in Theorem \ref{thm exist} and let $f$ satisfy Assumption \ref{assumption3} for $5R'$ and $L$.  For each $x \in B$ and $h \in [0,H')$ 
\begin{enumerate}
\item let $x' \in B_{R'}(x)$ and $\lambda \in \bbR$ be the unique solution to \eqref{eqn0} such that $|\lambda| \leq \tsfrac{11}{5} C_2 h$ (which exists by Theorem \ref{thm exist}), 
\item let $y \in B_{6R'}(x)$ be the unique solution to $y = x + h \tff(x,y,h)$ (which exists by \cite[Lem. 3.1]{nortonquispel1}), and 
\item let $x(\cdot)$ denote the exact solution to \eqref{p1} satisfying $x(t) = x$ for some $t \geq 0$.
\end{enumerate}
Also suppose that 
\begin{enumerate}
\setcounter{enumi}{3}
\item $\tff$ is such that the method $x \mapsto y$ defined by $y = x + h \tff(x,y,h)$ is of order $p$ for some $p \in \mathbb{N}$, i.e. when $f$ is sufficiently smooth there exist positive constants $C_3$ and $H_3 < H'$ such that  
\begin{equation}
\label{eqn30z}
	|y - x(t+h)| \leq C_3 h^{p+1} \qquad \mbox{for all $h \in [0,H_3]$ and all $x \in B$}.
\end{equation}
\end{enumerate}
Then the linear projection method defined by \eqref{eqn0} is also of order $p$, so that when $f$ is sufficiently smooth there exist positive constants $C_5$ and $H_5$ such that 
$$
	|x' - x(t+h)| \leq C_5 h^{p+1} \qquad \mbox{for all $h \in [0,H_5]$ and all $x \in B$}.
$$
\end{theorem}

\section{Existing linear projection methods}
\label{sec existing}

In this section we consider several linear projection methods that have been described and studied previously in the literature.  Our purpose is to show how all of these methods are special cases in our general framework, and hence our new theory also applies in these cases. 

We will need the following version of Banach's Fixed Point Theorem (also known as the Contraction Principle).   This version was also used in \cite{nortonquispel1} and is from \cite[Thm. 3.1.2 on p. 74]{fixedpointtheory}.

\begin{theorem}[Banach's Fixed Point Theorem]
\label{thm banach}
Let $(X,d)$ be a non-empty complete metric space.  Let $T : X \rightarrow X$ be a contraction on $X$, i.e. there exists a $q \in (0,1)$ such that 
$$
	d(T(x),T(y)) \leq q d(x,y) \qquad \mbox{for all $x,y \in X$}.
$$
Then there exists a unique fixed point $x^* \in X$ such that $T(x^*) = x^*$.  Furthermore, the fixed point can be found by iteration, $x^{n+1} = T(x^n)$ for $n=0,1,2,\dotsc$ with $x^0 \in X$ arbitrary.
\end{theorem}

\subsection{Example 1: (non-symmetric) standard projection method}   

In our notation, the \emph{(non-symmetric) standard projection method} described in \cite[Algorithm IV 4.2]{HLW} for $x \mapsto x'$ is defined by 
\begin{align*}
	y &= x + h \tgg(x,y,h) \\
	x' &= y + \lambda i(x') \mbox{ such that } I(x') = I(x),
\end{align*}
where $\tgg : \bbR^d \times \bbR^d \times [0,\infty) \rightarrow \bbR^d$ defines a map $x \mapsto y$ that is an arbitrary one-step method applied to \eqref{p1}.  Let $\Phi_h$ be the implicitly defined map so that $y = \Phi_h(x)$.  If we define $\tii : \bbR^d \times \bbR^d \times [0,\infty) \rightarrow \bbR^d$ and $\tff : \bbR^d \times \bbR^d \times [0,\infty) \rightarrow \bbR^d$ such that 
\begin{equation}
\label{eqn50}
	\tii(x,x',h) := i(x') \qquad \mbox{and} \qquad \tff(x,x',h) := \tgg(x,\Phi_h(x),h),
\end{equation}
for all $x,x' \in \bbR^d$ and $h \in [0,\infty)$ then the standard projection method is a linear projection method of the form \eqref{eqn0}.  

However, for computation the authors of \cite{HLW} suggest using \eqref{eqn0} with $\tii$ defined by 
\begin{equation}
\label{eqn51}
	\tii(x,x',h) := i(y) = i(\Phi_h(x))
\end{equation}
for all $x,x' \in \bbR^d$ and $h \in [0,\infty)$, instead of \eqref{eqn50} to make the system of equations easier to solve at each time step.  Strictly speaking, this method with $\tii$ given by \eqref{eqn51} instead of \eqref{eqn50} is a different projection method because the projection direction is different.  Let us refer to it as \emph{version 2} of the (non-symmetric) standard projection method.

To apply our new theory in Theorems \ref{thm exist} and \ref{thm order} we must determine what conditions $\tgg$ and $i$ must satisfy to ensure that $\tff$ and $\tii$ satisfy Assumptions \ref{assumption1} and \ref{assumption2}, respectively, for both versions of the standard projection method.  First, consider the first version of the standard projection method when $\tff$ and $\tii$ are defined by \eqref{eqn50}.  We must first prove the following lemma about existence, uniqueness and continuity of $\Phi_h$.

\begin{lemma}
\label{lem yexist}
For a bounded set $B \subset \bbR^d$, let $C_1$ be the constant from \eqref{eqn2} and suppose that $\tgg$ satisfies Assumption \ref{assumption1} for some positive constants $R_g$, $L_g$ and $H_g$.  

If $u \in B_{2R_g}(x)$ and $h < \min \{ H_g, \tsfrac{1}{6L_g}, \tsfrac{1}{4(C_1 + 1/6)R_g} \}$, then there exists a unique $y = \Phi_h(u) \in B_{R_g}(x)$ satisfying $y = u + h \tgg(u,y,h)$.  Moreover,
\begin{equation}
\label{eqn2s}
	|\Phi_h(u) - u| \leq \left( C_1 + \tsfrac{1}{6} + \tsfrac{3 L_g}{2 R_g} \right) h |i(x)|,
\end{equation}
and if $v \in B_{2R_g}(x)$ then
\begin{equation}
\label{eqn3s}
	|\Phi_h(u) - \Phi_h(v)| \leq \tsfrac{7}{5} |u-v|.
\end{equation}
\end{lemma}

\begin{proof}
Fix $x \in B$, $h < \min \lbrace H_g, \tsfrac{1}{6L_g}, \tsfrac{1}{4(C_1 + 1/6)R_g} \rbrace$ and $u,v \in B_{2R_g}(x)$.  We will apply Theorem \ref{thm banach} with $X := B_{R_g}(x)$ and $T(z):= u + h \tgg(u,z,h)$ for all $z \in X$.  To do so we must show that $T(z) \in X$ for any $z \in X$ and that $T$ is a contraction on $X$.  It is obvious that $X$ with the metric $|\cdot|$ (the usual Euclidean distance) is a non-empty complete metric space.  Let $z \in X$.  Then using Assumption \ref{assumption1} for $\tgg$, \eqref{eqn2}, $u \in B_{2R_g}(x) \subset B_{R_g}(x)$, $z \in B_{R_g}(x)$ and the bound on $h$ we get
\begin{equation}
\label{eqn1s}
\begin{split}
	|T(z) - x | &\leq |u - x| + h |\tgg(u,z,h)| \\
	&\leq |u-x| + h \left( |f(x)| + L_g |u-x| + L_g |z-x| + L_g h |i(x)| \right) \\
	&\leq \tsfrac{|i(x)|}{2R_g} + h \left( C_1 + \tsfrac{L_g}{2R_g} + \tsfrac{L_g}{R_g} + \tsfrac{1}{6} \right) |i(x)| \\
	&= \tsfrac{|i(x)|}{2R_g} + h \left( C_1 + \tsfrac{1}{6} \right) |i(x)| + h \tsfrac{3 L_g}{2R_g} |i(x)| 
	\leq \tsfrac{|i(x)|}{R_g}.
\end{split}
\end{equation}
Hence $T(z) \in X$.  For $z,z' \in X$, using Assumption \ref{assumption1} for $\tgg$ and $h \leq \tsfrac{1}{6L_g}$ we also get
$$
	|T(z) - T(z')| = h |\tgg(u,z,h) - \tgg(u,z',h)| \leq h L_g |z-z'| \leq \tsfrac{1}{6}|z-z'|
$$
so that $T$ is a contraction on $X$.  Therefore, by Theorem \ref{thm banach} there exists a unique $y \in B_{R_g}(x)$ satisfying $y = u + h \tgg(u,y,h)$.  Define $\Phi_h$ such that $\Phi_h(u) =y$.  

To get \eqref{eqn2s} we use a similar argument to \eqref{eqn1s}.  Finally, \eqref{eqn3s} follows from the following inequality where we have again used Assumption \ref{assumption1} for $\tgg$ and $h < \min \lbrace H_g, \tsfrac{1}{6L_g} \rbrace$
\begin{align*}
	|\Phi_h(u) - \Phi_h(v)| &\leq |u-v| + h | \tgg(u,\Phi_h(u),h) - \tgg(v,\Phi_h(v),h)| \\
	&\leq |u-v| + h \left( L_g |u - v| + L_g |\Phi_h(u) - \Phi_h(v)| \right) \\
	&\leq \tsfrac{7}{6} |u-v| + \tsfrac{1}{6} |\Phi_h(u) - \Phi_h(v)|.
\end{align*}
Hence $|\Phi_h(u) - \Phi_h(v)| \leq \tsfrac{7}{5} |u-v|$.
\end{proof}

Now we can prove that $\tff$ defined by \eqref{eqn50} satisfies Assumption \ref{assumption1} for some choice of $R$, $L$ and $H$.  

\begin{lemma}
\label{lem tff1}
For a bounded set $B \subset \bbR^d$, let $C_1$ be the constant from \eqref{eqn2} and suppose that $\tgg$ satisfies Assumption \ref{assumption1} for some positive constants $R_g$, $L_g$ and $H_g$.  Define $R := 2R_g$, 
$$
L := \max \left\{ \tsfrac{12}{5} L_g, \left(C_1 + \tsfrac{7}{6} + \tsfrac{3L_g}{2R_g} \right) L_g \right\}
	\quad \mbox{and} \quad
	H := \min \left\{ H_g, \tsfrac{1}{6L_g}, \tsfrac{1}{4(C_1 + 1/6) R_g} \right\}.
$$
Then $\tff : \bbR^d \times \bbR^d \times [0,\infty) \rightarrow \bbR^d$ defined by $\tff(x,x',h) := \tgg(x, \Phi_h(x), h)$ for all $x,x' \in \bbR^d$ and $h \in [0,\infty)$ satisfies Assumption \ref{assumption1} for $R$, $L$ and $H$.
\end{lemma}

\begin{proof}
Fix $x \in B$, $u,v,w \in B_R(x)$ and $h \in [0,H)$.  Using Assumption \ref{assumption1} for $\tgg$ and $\Phi_0(x) = x$ we get
$$
	\tff(x,x,0) = \tgg(x,\Phi_0(x),0) = \tgg(x,x,0) = f(x).
$$
Using Assumption \ref{assumption1} for $\tgg$ and Lemma \ref{lem yexist} (in particular \eqref{eqn3s}) we get 
\begin{align*}
	|\tff(u,v,h) - \tff(w,v,h)| &= |\tgg(u,\Phi_h(u),h) - \tgg(w,\Phi_h(w),h)| \\
	&\leq L_g |u-w| + L_g |\Phi_h(u) - \Phi_h(w)| \\
	&\leq L_g (1 + \tsfrac{7}{5}) |u-w| \\
	&\leq L |u-w|.
\end{align*}
Trivially, $|\tff(u,v,h) - \tff(u,w,h)| = 0$.  Finally, again using Assumption \ref{assumption1} for $\tgg$ and Lemma \ref{lem yexist} (in particular \eqref{eqn2s}) we get 
\begin{align*}
	|\tff(x,x,h) - \tff(x,x,0)| 
	&= |\tgg(x, \Phi_h(x),h) - \tgg(x,x,0) | \\
	&\leq L_g |\Phi_h(x) - x| + L_g h |i(x)| \\
	&\leq \left( L_g (C_1 + \tsfrac{1}{6} + \tsfrac{3L_g}{2R_g}) + L_g \right) h |i(x)| \\
	&\leq L h |i(x)|.
\end{align*}
\end{proof}

Now let us consider how we should choose $R$, $L$ and $H$ so that $\tii$ defined by \eqref{eqn50} or \eqref{eqn51} should satisfy Assumption \ref{assumption2}.  

\begin{lemma}
\label{lem tii1}
For a bounded set $B \subset \bbR^d$, let $R$ and $L$ be positive constants such that $i$ satisfies Assumption \ref{assumption3} for $R$ and $L$, and let $H$ be an arbitrary positive constant.  Then $\tii : \bbR^d \times \bbR^d \times [0,\infty) \rightarrow \bbR^d$ defined by $\tii(x,x',h) = i(x')$ for all $x,x' \in\bbR^d$ and $h \in [0,\infty)$ satisfies Assumption \ref{assumption2} for $R$, $L$ and $H$.
\end{lemma}

\begin{proof}
Since $i$ is locally Lipschitz, given arbitrary $R > 0$, $L$ exists.  The rest of the proof is trivial.
\end{proof}

\begin{lemma}
\label{lem tii1a}
For a bounded set $B \subset \bbR^d$, let $C_1$ be the constant from \eqref{eqn2}, and suppose that 
\begin{enumerate}
\item $\tgg$ satisfies Assumption \ref{assumption1} for some positive constants $R_g$, $L_g$ and $H_g$, and
\item $i$ satisfies Assumption \ref{assumption3} for $R_g$ and $L_i$ (given $R_g$, $L_i$ exists since $i$ is locally Lipschitz).
\end{enumerate}
Define $R := 2R_g$, 
\begin{align*}
	L &:= \max \left\{ \tsfrac{12}{5} L_g, \tsfrac{7}{5} L_i, \left(C_1 + \tsfrac{7}{6} + \tsfrac{3L_g}{2R_g} \right) \max \{ L_g, L_i \} \right\} \quad \mbox{and}\\
	H &:= \min \left\{ H_g, \tsfrac{1}{6L_g}, \tsfrac{1}{4(C_1 + 1/6) R_g} \right\}.
\end{align*}
Then $\tii: \bbR^d \times \bbR^d \times [0,\infty) \rightarrow \bbR^d$ defined by $\tii(x,x',h) := i(\Phi_h(x))$ for all $x,x' \in \bbR^d$ and $h \in [0,\infty)$ satisfies Assumption \ref{assumption2} for $R$, $L$ and $H$.  
\end{lemma}

\begin{proof}
Let $x \in B$, $u,v,w \in B_R(x)$ and $h \in [0,H)$.  Since $\Phi_0(x) = x$ it follows that $\tii(x,x,0) = i(\Phi_0(x)) = i(x)$.  Using Assumption \ref{assumption3} for $i$ and Lemma \ref{lem yexist} (in particular \eqref{eqn3s}) we get
$$
	|\tii(u,v,h) - \tii(w,v,h)| \leq L_i |\Phi_h(u) - \Phi_h(w)| \leq \tsfrac{7}{5} L_i |u-w| \leq L |u-w|.
$$
Trivially, we have $|\tii(u,v,h) - \tii(u,w,h)| = 0$, and using $x = \Phi_0(x)$, Assumption \ref{assumption3} for $i$ and Lemma \ref{lem yexist} (in particular \eqref{eqn2s}) we get
$$
	|\tii(x,x,h) - \tii(x,x,0)| \leq L_i |\Phi_h(x) - x| \leq L_i \left( C_1 + \tsfrac{1}{6} + \tsfrac{3L_g}{2R_g} \right) h |i(x)| \leq L h |i(x)|.
$$
\end{proof}
\subsection{Example 2:  symmetric projection method}  
\label{sec ex2}

It is perhaps surprising that the symmetric projection method from \cite[\S V.4.1]{HLW} (originally in \cite{hairer00}) may also be written in the form of \eqref{eqn0}.  In our notation, the symmetric projection method described in \cite[\S V 4.1]{HLW} for $x \mapsto x'$ is defined by: Given $x \in \bbR^d$ and $h \in [0,\infty)$, compute $y,z,x' \in \bbR^d$ and $\mu \in \bbR$ such that 
\begin{equation}
\label{eqn4s}
\begin{split}
	y &= x + \mu i(x), \\
	z &= y + h \tgg(y,z,h), \\
	x' &= z + \mu i(x') \mbox{ and } I(x') = I(x), 
\end{split}
\end{equation}
where $\tgg : \bbR^d \times \bbR^d \times [0,\infty) \rightarrow \bbR^d$ is such that $y \mapsto z$ defined by $z = y + h \tgg(y,z,h)$ is any symmetric one-step method applied to \eqref{p1}.  If we let $\lambda = 2 \mu$ and eliminate $y$ and $z$ from \eqref{eqn4s} then we get: Given $x \in \bbR^d$ and $h \in [0,\infty)$ compute $x' \in \bbR^d$ and $\lambda \in \bbR$ such that 
$$
	x' = x + h \tgg\left( x + \tsfrac{\lambda}{2} i(x), x' - \tsfrac{\lambda}{2}i(x') ,h\right) + \lambda \tsfrac{ i(x) + i(x') }{2} \qquad \mbox{and} \qquad I(x') = I(x).
$$
If we let $\Psi$ be the implicitly defined mapping so that $\lambda = \Psi(x,x',h)$ where $\lambda$ satisfies
$$
	\lambda = - h \frac{ \tgg\left(x + \tsfrac{\lambda}{2} i(x), x' - \tsfrac{\lambda}{2} i(x'),h\right) \cdot \bii(x,x')}{\tsfrac{i(x) + i(x')}{2} \cdot \bii(x,x')}
$$
where $\bii$ is an arbitrarily chosen discrete gradient of $I$, then we see that the symmetric projection method is equivalent to \eqref{eqn0} if we define $\tii : \bbR^d \times \bbR^d \times [0,\infty) \rightarrow \bbR^d$ and $\tff : \bbR^d \times \bbR^d \times [0,\infty) \rightarrow \bbR^d$ by 
\begin{equation}
\label{eqn53}
\begin{split}
	\tii(x,x',h) &:= \frac{i(x) + i(x')}{2}, \\
	\tff(x,x',h) &:= \tgg\left(x + \tsfrac{\Psi(x,x',h)}{2} i(x), x' - \tsfrac{\Psi(x,x',h)}{2} i(x'),h\right)
\end{split}
\end{equation}
for all $x,x'\in \bbR^d$ and $h \in [0,\infty)$.

It turns out that the fact that $\tgg$ satisfies Assumption \ref{assumption1} is sufficient to ensure that $\tff$ and $\tii$ satisfy Assumptions \ref{assumption1} and \ref{assumption2} respectively.  This will ensure that we are able to apply Theorems \ref{thm exist} and \ref{thm order} to the symmetric projection method.  However, verifying that this is true is quite technical, so we have included it only as an appendix.

\subsection{Example 3: Methods of Dahlby, Owren and Yaguchi}  
\label{sec ex3}

In \cite{dahlby11} Dahlby et al. describe two projection methods.  In our notation, given an arbitrary discrete gradient $\bii$ of $I$, then the first of their methods (see \cite[eq. 2.1]{dahlby11}) is defined by
\begin{equation}
\label{eqn8}
	y = \Phi_h(x), \qquad x' = x + P (y-x)	
\end{equation}
where $\Phi_h : \bbR^d \rightarrow \bbR^d$ defines an arbitrary one-step method for solving \eqref{p1}, and $P$ is a projection onto $\opspan \{ \bii(x,x') \}^\perp$.  In \cite[\S 2.2]{dahlby11} the projection matrix $P$ is defined as
\begin{equation}
\label{eqn8a}
	P = I - \frac{\bii(x,x') \bii(x,x')^T}{\bii(x,x')^T \bii(x,x')},
\end{equation}
which is the orthogonal projection matrix onto $\opspan \{ \bii(x,x') \}^\perp$.
If we define $\tff : \bbR^d \times \bbR^d \times [0,\infty) \rightarrow \bbR^d$ and $\tii : \bbR^d \times \bbR^d \times [0,\infty) \rightarrow \bbR^d$ by
\begin{equation}
\label{eqn57}
	\tff(x,x',h) := \frac{\Phi_h(x) - x}{h} 
	\qquad \mbox{and} \qquad
	\tii(x,x',h) := \bii(x,x')
\end{equation}
for all $x,x' \in \bbR^d$ and $h \in [0,\infty)$ then it is easy to see that \eqref{eqn8} is the same method as \eqref{eqn2nd}, so it is a special case of our general linear projection method.  

The second of the projection methods by Dahlby et al.  (see \cite[eq. 2.2]{dahlby11}) is, in our notation and given an arbitrary discrete gradient $\bii$ of $I$, defined by
\begin{equation}
\label{eqn57a}
	x' = x + h P \tgg(x,x',h)
\end{equation}
where $P$ is the same projection matrix as above and $\tgg : \bbR^d \times \bbR^d \times [0,\infty) \rightarrow \bbR^d$ is such that the map $x \mapsto y$ defined by $y = x + h \tgg(x,y,h)$ is an arbitrary one-step method for solving \eqref{p1}.  By defining $\tff$ and $\tii$ such that 
\begin{equation}
\label{eqn58}
	\tff(x,x',h) = \tgg(x,x',h) 
	\qquad \mbox{and} \qquad
	\tii(x,x',h) = \bii(x,x')
\end{equation}
for all $x,x' \in \bbR^d$ and $h \in [0,\infty)$ we see that \eqref{eqn57a} is the same as \eqref{eqn2nd}, so it is another special case of our general linear projection method.

Dahlby et al. call their methods `discrete gradient methods' because the methods are constructed using a discrete gradient.  We think it is more appropriate to describe these methods as projection methods.  However, our theory (Theorem \ref{thm equiv1}) has established that they may also be expressed in the form of \eqref{disc1}, for which we use the term `discrete gradient method'.

It is a relatively simple task to show that if for a given bounded set $B \subset \bbR^d$, $\tgg$ satisfies Assumption 1 for some positive constants $R_g$, $L_g$ and $H_g$, and $\bii$ is a discrete gradient of $I$ satisfying Assumption \ref{assumption2} for some positive constants $R_{\bii}$ and $L_{\bii}$, then $\tff$ and $\tii$ defined by \eqref{eqn58} satisfy Assumptions \ref{assumption1} and \ref{assumption2}, respectively, for some positive constants $R$, $L$ and $H$.  For this reason and for the sake of brevity we omit the details.  In the case when $\tff$ and $\tii$ are defined by \eqref{eqn57} we must make suitable assumptions about the method defined by $\Phi_h$ for a similar result to hold.




\section{Equivalence in the multiple first integrals case}
\label{sec equiv2}

Now let us consider the case when \eqref{p1} has multiple preserved integrals.  Suppose that \eqref{p1} preserves $M$ first integrals $I_1, \dotsc, I_M$, i.e. there exist $I_m : \bbR^d \rightarrow \bbR$ for $m=1,\dotsc,M$ such that for all $t \geq 0$,
$$
	x(t) \in \mathcal{M}_{x_0} := \{ z \in \bbR^d : I_m(z) = I_m(x_0) \mbox{ for all } m = 1,\dotsc,M \}.
$$
For each $m$ we use the notation, $i_m := \nabla I_m$.
Recall from Section \ref{sec intro} the algorithm for computing one step $x \mapsto x'$ of a projection method: Given $x \in \bbR^d$ and $h \in [0,\infty)$
\begin{enumerate}
\item compute $y$ such that $y = x + h \tff(x,y,h)$, 
\item compute $x' \in \mathcal{M}_x$ by projecting $y$ onto $\mathcal{M}_x$.
\end{enumerate}
For the general linear projection case when we have multiple integrals to preserve we first define $M$ directions for the projection, i.e. $\tii_m : \bbR^d \times \bbR^d \times [0,\infty) \rightarrow \bbR^d$ for $m = 1, \dotsc, M$, and then replace step 2 with
\begin{enumerate}
\setcounter{enumi}{1}
\item compute $x' \in \bbR^d$ by solving $x'  = y + A \lambda$ such that $x' \in \mathcal{M}_x$ for $x' \in \bbR^d$ and $\lambda \in \bbR^M$, where $A = [ \tii_1 \dotsi \tii_M] \in \bbR^{d \times M}$ and $\tii_m = \tii_m(x,x',h)$.
\end{enumerate}

We say that this type of projection is a linear projection because $x'-y$ is a linear combination of the projection directions $\tii_m$, i.e. $(x'-y) \in \opspan \{ \tii_1, \dotsc, \tii_M \}$.  As in Section \ref{sec intro} for the single integral case, we can write the two step general linear projection method algorithm in one line by eliminating $y$ and generalising $\tff$.  The method is: Given $x \in \bbR^d$ and $h \in [0,\infty)$ compute $x' \in \bbR^d$ and $\lambda \in \bbR^M$ such that
\begin{equation}
\label{eqn0m}
	x' = x + h \tff(x,x',h) + A \lambda \qquad \mbox{and} \qquad I_m(x') = I_m(x) \mbox{ for all $m$}.
\end{equation}

\subsection{Oblique Projections}

Before we present our theorem showing the equivalence between linear projection methods and discrete gradient methods for ODEs with multiple integrals we need to introduce \emph{oblique projection matrices}.

The type of projection described in detail in linear algebra textbooks is usually orthogonal projection, e.g. \cite[Lecture 6]{trefethenbau} and \cite[\S 5.13]{meyer}.  For a set of linearly independent vectors $a_1, \dotsc, a_M \in \bbR^d$ (not necessarily orthogonal), the orthogonal projection matrix that maps $\bbR^d$ onto the subspace $\cA := \opspan \lbrace a_1,\dotsc,a_M \rbrace$ along the subspace $\cA^\perp$ is given by 
$$
	Q = A (A^T A)^{-1} A^T, \qquad \mbox{where $A = [ a_1 \dotsi a_M] \in \bbR^{d \times M}$}.
$$
Since $x = Q x + (I-Q) x$ is the unique decomposition of $x \in \bbR^d$ into vectors in $\cA$ and $\cA^\perp$ (see e.g. \cite[p. 43]{trefethenbau}), it follows that 
$$
	Q_\perp = I - Q = I - A(A^TA)^{-1} A^T
$$
is the projection matrix onto $\cA^\perp$ along $\cA$.  

We would like to consider the generalisation of these two projections where the range of the projection is decoupled from the direction of the projection.  For this we define an \emph{oblique projection} (see e.g. \cite{meyer}).  Note that the space along which a projection projects is the null space of the projection and we may define an oblique projection by specifying its range and null space.  For a projection $R$ let $a_1,\dotsc,a_M$ be a basis for the range of $R$, and let $b_1,\dotsc,b_M$ be a basis for the orthogonal complement of the null space of $R$, so that
\begin{align*}
	\oprange(R) &= \cA, \\
	\opnull(R) &= \cB^\perp,
\end{align*}
where $\cA := \opspan  \{ a_1, \dotsc, a_M \}$ and $\cB := \opspan \{ b_1, \dotsc, b_M \}$.  Then the oblique projection matrix $R$ is given by the formula \cite[eq. (7.10.39) on p. 634]{meyer},
$$
	R = A (B^T A)^{-1} B^T
$$
where $A = [a_1 \dotsi a_M] \in \bbR^{d \times M}$ and $B = [b_1 \dotsi b_M] \in \bbR^{d \times M}$.  To ensure that $B^T A$ is invertible and $R$ exists we also need that $\cA$ and $\cB^\perp$ are complementary subspaces of $\bbR^d$ (i.e. $\cA + \cB^\perp = \bbR^d$ and $\cA \cap \cB^\perp = \{ 0 \}$, see e.g. \cite[p. 383]{meyer}).  More generally, in the terminology of \cite{meyer}, if $M = A B^T$ is a full rank factorisation of $M \in \bbR^{d\times M}$ and if $\oprange(M)$ and $\opnull(M)$ are complimentary subspaces, then $R = A(B^TA)^{-1}B^T$ is the projection onto $\oprange(M)=\oprange(A)$ along $\opnull(M) = \opnull(B^T)$.  Note that $A$ and $B^T$ have full rank if the columns of $A$ are linearly independent, and the rows of $B^T$ are linearly independent (see e.g. \cite[p. 218]{meyer}).  

The following proposition will help us decide whether or not $\cA$ and $\cB^\perp$ are complementary subspaces.

\begin{proposition}
\label{prop0}
With $a_1,\dotsc,a_M$, $b_1,\dotsc,b_M$, $\cA$, $\cB$, $A$ and $B$ defined as above, then
$\cA$ and $\cB^\perp$ are complementary subspaces if and only if $B^T A$ is invertible.	
\end{proposition}

\begin{proof}
From \cite[p. 383]{meyer} we have that $\cA$ and $\cB^\perp$ are complimentary subspaces if and only if for any $x \in \bbR^d$ there exists a unique decomposition $x = x_{\cA} + x_{\cB^\perp}$ where $x_{\cA} \in \cA$ and $x_{\cB^\perp} \in \cB^\perp$.

Assume that $\cA$ and $\cB^\perp$ are complementary subspaces.  Then there exists a unique decomposition $x = x_{\cA} + x_{\cB^\perp}$ and since $ a_1,\dotsc,a_M $ is a basis for $\cA$ there exists a unique $v \in \bbR^M$ such that $x_{\cA} = A v$.  Therefore, $B^T x = B^T x_{\cA} = B^T A v$ and since $v$ is uniquely determined given $x$, $B^T A$ is invertible.    

Conversely, suppose $B^T A$ is invertible.  Then the matrix $R = A(B^T A)^{-1} B^T$ is well defined and for a given $x \in \bbR^d$, $x_{\cA} := R x \in \cA$ and $x_{\cB^\perp} := (I - R)x \in \cB^\perp$ defines a decomposition $x = x_{\cA} + x_{\cB^\perp}$.  To complete the proof we must show that this decomposition is unique.  Suppose $x = y_{\cA} + y_{\cB^\perp}$ where $y_{\cA} \in \cA$ and $y_{\cB^\perp} \in \cB^\perp$ defines another decomposition of $x$.  There exists a unique $w \in \bbR^M$ such that $y_{\cA} = A w$.  Then $B^T x = B^T y_{\cA} = B^T A w$ and hence $w = (B^T A)^{-1} B^T x$.  Substituting this into $y_{\cA} = A w$ we get $y_{\cA} = R x = x_{\cA}$ and $y_{\cB^\perp} = x_{\cB^\perp}$ and the decomposition is unique.
\end{proof}

An obvious choice for $B$ so that $B^TA$ is invertible is $B=A$.  Since $A$ has full rank, $A^T A$ is positive definite and invertible.  But this corresponds to orthogonal projection.  More generally, if $B$ is sufficiently ``close'' to $A$ then $B^T A$ is positive definite and hence invertible.  If for any $m$, $b_m \in \cA^\perp$ then $B^T A$ is not invertible since it has a column with all zeros.  

Another projection matrix $R_\bot$ with
\begin{align*}
	\oprange(R_\bot) &= \cB^\perp\\
	\opnull(R_\bot) &= \cA
\end{align*}
may be defined by
$$
	R_\perp := I - R = I - A (B^T A)^{-1} B^T.
$$
We will use $R_\bot$ to define an alternative formulation for linear projection methods for ODEs with multiple first integrals.

\subsection{Equivalent formulation using an oblique projection matrix}

A general method for solving \eqref{p1},
\begin{equation}
\label{eqna1}
	x' = x + h \tff(x,x',h)
\end{equation}
is integral preserving if and only if $\tff \in \opspan \{ \bii_1, \dotsc, \bii_M \}^\perp$ where each $\bii_m$ is a discrete gradient of $I_m$.  This fact follows from the definition of a discrete gradient: For each $m$, if $\bii_m$ is a discrete gradient of $I_m$, then $\tff \perp \bii_m$ if and only if 
\begin{equation}
\label{eqna2}
	I_m(x') - I_m(x) = (x' - x) \cdot \bii_m = h \tff \cdot \bii_m = 0.
\end{equation}
However, in general, we do not have $\tff \in \opspan \{ \bii_1, \dotsc, \bii_M \}^\perp$.  Therefore, a way of constructing an integral preserving matrix is to modify \eqref{eqna1} to
$$
	x' = x + h P \tff(x,x',h)
$$
where $P = P(x,x',h)$ is chosen to be a projection matrix with $\oprange(P) = \opspan \{ \bii_1, \dotsc, \bii_M \}^\perp$ so that $P \tff \in \opspan \{ \bii_1, \dotsc, \bii_M \}^\perp$.  It turns out that constructing an integral preserving method in this way is equivalent to a general linear projection method of the form \eqref{eqn0m}.  This equivalence is formalised in the following theorem and is an extension to our earlier Theorem \ref{thm equiv1} (in particular the equivalence between \eqref{eqn1st} and \eqref{eqn2nd}).

\begin{theorem}
\label{thm equiv1m}
Let $\tff : \bbR^d \times \bbR^d \times [0,\infty) \rightarrow \bbR^d$ be a consistent approximation of $f$ and for each $m = 1, \dotsc, M$ let $\tii_m : \bbR^d \times \bbR^d \times [0,\infty) \rightarrow \bbR^d$ be a consistent approximation of $i_m$ and let $\bii_m$ be an arbitrary discrete gradient of $I_m$.  Define 
\begin{equation}
\label{eqna7}
	P := I - A (B^T A)^{-1} B^T
\end{equation}
where
$$
	A := [\tii_1 \dotsi \tii_M] \in \bbR^{d \times M}\qquad \mbox{and} \qquad 
	B := [\bii_1 \dotsi \bii_M] \in \bbR^{d \times M},
$$
and $\tii_m = \tii_m(x,x',h)$ and $\bii_m = \bii_m(x,x')$.  Assume that 
\begin{enumerate}
\item each of the two methods below have uniquely defined $x' \in \bbR^d$ and $\lambda \in \bbR^M$ for sufficiently small $h$, 
\item $\{ \tii_1,\dotsc,\tii_M \}$ and $\{ \bii_1, \dotsc, \bii_M \}$ are linearly independent sets (so that $A$ and $B^T$ have full rank), and 
\item $\tilde{\cS}$ and $\bar{\cS}^\perp$ are complementary subspaces of $\bbR^d$, where $\tilde{\cS} = \opspan \{ \tii_1,\dotsc, \tii_M \}$ and $\bar{\cS} = \opspan \{ \bii_1, \dotsc, \bii_M \}$ (i.e. $\tilde{\cS} + \bar{\cS}^\perp = \bbR^d$ and $\tilde{\cS} \cap \bar{\cS}^\perp = \{ 0 \}$).
\end{enumerate}
Then the following expressions describe the same linear projection method.
\begin{equation}
\label{eqna1st}
	x' = x + h \tff(x,x',h) + A \lambda \mbox{ such that $I_m(x') = I_m(x)$ for all $m=1,\dotsc,M$,}
\end{equation}
and
\begin{equation}
\label{eqna2nd}	
	x' = x + h P \tff(x,x',h).
\end{equation}
\end{theorem}

\begin{proof}
Conditions 2 and 3 in the theorem ensure that $P$ exists ($B^T A$ is invertible), $\oprange(P) = \bar{\cS}^\perp$ and $\opnull(P) = \tilde{\cS}$ (see discussion in previous section about oblique projection matrices).

For given $x \in \bbR^d$ and $h$ sufficiently small, suppose $x' \in \bbR^d$ and $\lambda \in \bbR^M$ satisfy \eqref{eqna1st}.  For each $m$, since $\bii_m$ is a discrete gradient,
$$
	0 = I_m(x') - I_m(x) = \bii_m(x,x')^T (x'-x) = h \bii_m^T \tff + \bii_m^T A \lambda.
$$
Therefore,
$$
	0 = h B^T \tff + B^T A \lambda \qquad \Rightarrow \qquad
	\lambda = - h (B^T A)^{-1} B^T \tff.
$$
Substituting this into \eqref{eqna1st} we get
$$
	x' = x + h \tff - h A (B^T A)^{-1} B^T \tff = x + h P \tff,
$$
so $x'$ satisfies \eqref{eqna2nd}.

Conversely, for given $x \in \bbR^d$ and $h$ sufficiently small, suppose $x' \in \bbR^d$ satisfies \eqref{eqna2nd}.  We know that $\oprange (P) = \bar{\cS}^\perp$.  Therefore, $P \tff \in \bar{\cS}^\perp$ and by \eqref{eqna2} we have $I_m(x') = I_m(x)$ for all $m = 1,\dotsc,M$.  Finally, if we define $\lambda = - h (B^T A)^{-1} B^T \tff$ then
$$
	x' = x + h P \tff = x + h \tff - h A (B^T A)^{-1} B^T \tff = x + h \tff + A \lambda
$$
and $x'$ and $\lambda$ satisfy \eqref{eqna1st} as required.
\end{proof}

Since $\tii_m$ and $\bii_m$ are consistent approximations of $i_m$ it follows that if $i_m(x)$ for $m = 1,\dotsc,M$ are linearly independent then for sufficiently small $h$, both $\{ \tii_1,\dotsc,\tii_M \}$ and $\{\bii_1, \dotsc, \bii_M \}$ are linearly independent sets.  Moreover, for small $h$ the matrix $B$ is ``close'' to $A$ and the property that $\tilde{\cS}$ and $\bar{\cS}^\perp$ are complementary subspaces of $\bbR^d$ is satisfied (see discussion after Proposition \ref{prop0}).

Now that we have established an alternative formulation for general linear projection methods we can explore their relationship to discrete gradient methods.




\subsection{Equivalence with discrete gradient methods}

In this section let us consider discrete gradient methods for preserving more than one integral.  Our aim is to construct a general discrete gradient method to approximate the solution to \eqref{p1} such that $M$ integrals are simultaneously preserved, and then determine which discrete gradient methods are equivalent to linear projection methods.  

According to \cite[Prop. 2.14]{MQR99} (see also \cite{quispeldyt} for the two integral case), we may write \eqref{p1} as 
\begin{equation}
\label{eqn60}
	\dot{x}_j = S_{j j_1 j_2 \; \dotsi j_M} i^1_{j_1} i^2_{j_2} \dotsi i^M_{j_M} \qquad \mbox{for each $j=1,\dotsc,d$},
\end{equation}
using Einstein's summation principle for repeated indices, where $i^m := \nabla I_m$ for each $m$ and
$$
	S = \frac{1}{\det (D^T D)} f \wedge i^1 \wedge \dotsi \wedge i^M,
$$
where $D: \bbR^d \rightarrow \bbR^{d \times M}$ is defined by $D = [i^1 \dotsi i^M]$.  See e.g. \cite[Chap. 1]{darling} or \cite{greub} for a definition of the anti-symmetric $\wedge$ product from exterior algebra.

Based on the expression for the ODE given in \eqref{eqn60} we can write down a general discrete gradient method for solving \eqref{p1}.  Let $\wS = \wS(x,x',h)$ be an anti-symmetric consistent approximation of $S$ and for each $m$ let $\bii^m = \bii^m(x,x')$ be a discrete gradient of $I_m$.  Then, the method $x \mapsto x'$ is defined as: Given $x \in \bbR^d$ and $h \in [0,\infty)$, $x' \in \bbR^d$ satisfies 
\begin{equation}
\label{eqn61}
	\frac{x_j'-x_j}{h} = \wS_{j j_1 j_2 \; \dotsi j_M} \bii^1_{j_1} \dotsi \bii^M_{j_M} \qquad \mbox{for each $j=1,\dotsc,d$}.
\end{equation}

For particular choices of $\wS$, this discrete gradient method is equivalent to a linear projection method.  To prove our result we will need the following proposition (see \cite[eq. 5.10 on p. 106]{greub} where the pair of dual spaces are $\bbR^d$ and itself with the usual Euclidean inner product).

\begin{proposition}
\label{prop1}
For arbitrary $M,d \in \mathbb{N}$, let $U,V \in \bbR^{d \times M}$ be two matrices with columns $u^1,\dotsc,u^M \in \bbR^d$ and $v^1,\dotsc,v^M \in \bbR^d$ respectively.  Then 
$$
	(u^1 \wedge u^2 \wedge \dotsi \wedge u^M)_{j_1 j_2 \;\dotsi j_M} v^1_{j_1} \dotsi v^M_{j_M} = \det (V^T U).
$$
\end{proposition}

\begin{theorem}
\label{thm equiv3}
Let $\tff = \tff(x,x',h)$ be a consistent approximation of $f(x)$, and for each $m=1,\dotsc,M$ let $\tii^m = \tii^m(x,x',h)$ be a consistent approximation of $i^m(x)$ and let $\bii^m = \bii^m(x,x')$ be a discrete gradient of $I_m(x)$.  Define 
$$
	A = [ \tii^1 \dotsi \tii^M] \in \bbR^{d \times M} \qquad\mbox{and}\qquad B = [ \bii^1 \dotsi \bii^M] \in \bbR^{d \times M}.
$$
Assume that
\begin{enumerate}
\item the discrete gradient method and the linear projection method defined below have uniquely defined $x' \in \bbR^d$ and $\lambda \in \bbR^M$ for sufficiently small $h$,
\item $\{\tii^1,\dotsc,\tii^M \}$ and $\{ \bii^1,\dotsc,\bii^M \}$ are linearly independent sets (so that $A$ and $B$ have full rank), and
\item $\tilde{\cS}$ and $\bar{\cS}^\perp$ are complementary subspaces of $\bbR^d$ where $\tilde{\cS} = \opspan \{ \tii^1,\dotsc, \tii^M \}$ and $\bar{\cS} = \opspan \{ \bii^1, \dotsc, \bii^M \}$.
\end{enumerate}
If we define
\begin{equation}
\label{eqn62}
	\wS = \frac{1}{\det(B^T A)} \tff \wedge \tii^1 \wedge \dotsi \wedge \tii^M,
\end{equation}
then the discrete gradient method defined by \eqref{eqn61} and \eqref{eqn62} is equivalent to the linear projection method defined by \eqref{eqna1st} or \eqref{eqna2nd}.
\end{theorem}

\begin{proof}
To show that these methods are the same we must show that 
$$
	(P \tff)_j = \wS_{j j_1 j_2 \; \dotsi j_M} \bii^1_{j_1} \bii^2_{j_2} \dotsi \bii^M_{j_M}
$$
for each $j = 1,\dotsc,M$, where $P = I - A(B^T A)^{-1} B^T$.  Equivalently, we can show that 
$$
	(P \tff) \cdot v = \wS_{j j_1 j_2 \; \dotsi j_M} v_j \bii^1_{j_1} \bii^2_{j_2} \dotsi \bii^M_{j_M},
$$
for any $v \in \bbR^d$.  Let $v$ be an arbitrary vector in $\bbR^d$.  Using Proposition \ref{prop1} and expanding the determinant along the first row we get
\begin{align*}
	(\tff \wedge \tii^1 \wedge \dotsi \wedge \tii^M)_{j j_1 j_2 \; \dotsi j_M} v_j \bii^1_{j_1} \bii^2_{j_2} \dotsi \bii^M_{j_M} 
	&= \det \left( [ v \; \bii^1 \dotsi \bii^M]^T [ \tff \; \tii^1 \dotsi \tii^M] \right) \\
	&= \det \left( \left[ \begin{array}{c|ccc}
		v \cdot \tff & v \cdot \tii^1 & \dotsi & v \cdot \tii^M \\ \hline
		B^T \tff & & B^T A & 
		\end{array} \right] \right) \\
	&= (v \cdot \tff) \det( B^T A) + \sum_{j=1}^M (-1)^j (v \cdot \tii^j) \det( B^T \widetilde{A}_j ),
\end{align*}
where $\widetilde{A}_j = [\tff \; \tii^1 \dotsi \tii^{j-1} \tii^{j+1} \dotsi \tii^M ] \in \bbR^{d \times M}$.  Using the fact that the determinant of a matrix is anti-symmetric (each column swap introduces a factor of $-1$) it follows that
$$
	\det(B^T \widetilde{A}_j ) = (-1)^j \det(B^T A_j)
$$
where $A_j = [\tii^1 \dotsi \tii^{j-1} \tff \; \tii^{j+1} \dotsi \tii^M] \in \bbR^{d \times M}$, i.e. the matrix $A$ with $\tii^j$ replaced by $\tff$.  Hence, using the two identities above and Cramer's Rule (see e.g. \cite[p. 476]{meyer}) we get 
\begin{align*}
	\wS_{i j_1 j_2 \; \dotsi j_M} v_i \bii^1_{j_1} \bii^2_{j_2} \dotsi \bii^M_{j_M} 
	&= \frac{1}{\det (B^T A)} \left( (v \cdot \tff) \det( B^T A) - \sum_{j=1}^M (v \cdot \tii^j) \det( B^T A_j )\right)  \\
	&= \left( \tff - \sum_{j=1}^M \frac{\det(B^T A_j)}{\det(B^T A)} \tii^j \right) \cdot v \\
	&= \left( \tff - \sum_{j=1}^M \left[ (B^T A)^{-1} B^T \tff \right]_j \tii^j \right) \cdot v \\
	&= \left( \tff - A (B^T A)^{-1} B^T \tff \right) \cdot v = (P \tff) \cdot v.
\end{align*}
\end{proof}

If we restrict ourselves to the situation where only two integrals are preserved, $I$ and $J$, with $i := \nabla I$ and $j := \nabla J$, then \eqref{p1} may be written as (see \cite{quispeldyt})
\begin{equation}
\label{eqna3}
	\dot{x}_l = \sum_{m,n=1}^d S_{l m n} i_m j_n \qquad \mbox{for each $l=1,\dotsc,d$},
\end{equation}
where $S_{lmn}$ is an anti-symmetric tensor given by
$$
	S_{lmn} = \frac{1}{N} \left| \begin{array}{ccc}
		f_l & i_l & j_l \\
		f_m & i_m & j_m \\
		f_n & i_n & j_n \end{array} \right| 
	\qquad
	\mbox{and}
	\qquad	
	N = (i\cdot i)(j\cdot j) - (i\cdot j)^2.
$$
The general discrete gradient methods for \eqref{eqna3} that will preserve both $I$ and $J$ are
\begin{equation}
\label{eqna4}
	\frac{x_l' - x_l}{h} = \sum_{m,n=1}^d \wS_{lmn} \bii_m \bjj_n
\end{equation}
where $\bii$ and $\bjj$ are discrete gradients of $I$ and $J$ respectively, and $\wS$ is a skew-symmetric consistent approximation of $S$.  If we define 
$$
	\wS_{lmn} = \frac{1}{\widetilde{N}} \left| \begin{array}{ccc}
		\tff_l & \tii_l & \tjj_l \\
		\tff_m & \tii_m & \tjj_m \\
		\tff_n & \tii_n & \tjj_n \end{array} \right| 
	\qquad
	\mbox{and}
	\qquad	
	\tilde{N} =(\tii \cdot \bii) (\tjj \cdot \bjj) - (\tii \cdot \bjj)(\bii \cdot \tjj)
$$
where $\tff$ is a consistent approximation of $f$, $\tii$ and $\tjj$ are consistent approximations of $i$ and $j$ respectively, then this discrete gradient method is a linear projection method.

We remark that \eqref{eqna1st} and \eqref{eqna2nd} do not depend on any discrete gradients of $I_m$.  Therefore, we may conclude from Theorem \ref{thm equiv3} that each projection method (defined by the choice of projection directions) is equivalent to a class of discrete gradient methods where the approximate solution values at each time step are independent of the particular choices of discrete gradients used in the discrete gradient methods.

\subsection{Existence, uniqueness and order of accuracy}

For the single preserved integral case we could use theory from discrete gradient methods to prove, under certain local Lipschitz continuity and consistency conditions, the existence of a unique solution $x'$ at each time step of a projection method for sufficiently small $h$.  We were also able to show under the same conditions that a projection method retained the same order of accuracy as the underlying method.  For the multiple integral case we cannot do this because these results for discrete gradient methods are not yet available.  We do not anticipate that extending these results to the multiple integral case poses any real difficulty, except that a proof may be very lengthy.  In Section \ref{sec numerics} we will try to test numerically whether or not it is correct to assume that these results hold in the multiple integral case.

\subsection{Special cases of projection methods}

As in Section \ref{sec existing} we can show how our expression for a general linear projection method \eqref{eqn0m} for preserving multiple first integrals encompasses all existing (as far as we are aware) projection methods.  Unlike Section \ref{sec existing} however, we will not go into all of the technicalities regarding existence of a unique solution for sufficiently small time step.  Once again, the trick to seeing how other projection methods fit into our framework is to make the right choice for the projection directions $\tii_m$.  

{\noindent\bf Example 1 revisited: (non-symmetric) standard projection method.}

The first method we consider is again the (non-symmetric) standard projection method described in \cite[Algorithm IV 4.2]{HLW}.  In our notation, their method in the multiple preserved integral case for $x \mapsto x'$ is defined by solving the following system of equations for $x' \in \bbR^d$ and $\lambda \in \bbR^M$, given $x \in \bbR^d$ and $h \in [0,\infty)$,
\begin{align*}
	y &= x + h \tgg(x,y,h), \\
	x' &= y + A \lambda \qquad \mbox{and $I_m(x') = I_m(x)$ for every $m = 1, \dotsc,M$},
\end{align*}
where the map $x \mapsto y$ defined by $y = x + h \tgg(x,y,h)$ defines an arbitrary one-step method applied to \eqref{p1} and $A = [i_1(x') \dotsi i_M(x')] \in \bbR^{d \times M}$.  If we let $\Phi_h$ be the implicitly defined map so that $y = \Phi_h(x)$ then this method has the form of \eqref{eqn0m} if for each $m$ we define $\tii_m : \bbR^d \times \bbR^d \times [0,\infty) \rightarrow \bbR^d$ and $\tff : \bbR^d \times \bbR^d \times [0,\infty) \rightarrow \bbR^d$ such that
$$
	\tii_m(x,x',h) := i_m(x') \qquad \mbox{and} \qquad \tff(x,x',h) := \tgg(x,\Phi_h(x),h)
$$
for all $x,x' \in \bbR^d$ and $h \in [0,\infty)$.  In \cite{HLW}, the authors suggest instead using $\tii_m(x,x',h) := i_m(y)$ to reduce the number of evaluations of $i_m(\cdot)$ required when solving the system of equations at each step using a simplified Newton method.

{\noindent\bf Example 2 revisited: symmetric projection method.}

The multiple first integral version of the symmetric projection method (see \cite[\S V.4.1]{HLW} or \cite{hairer00}), in our notation for $x \mapsto x'$, is: Given $x \in \bbR^d$ and $h \in [0,\infty)$, compute $y,z,x' \in \bbR^d$ and $\mu \in \bbR^M$ such that
\begin{align*}
	y &= x + A'' \mu, \\
	z&= y + h \tgg(y,z,h), \\
	x'&= z + A' \mu \qquad \mbox{and $I_m(x') = I_m(x)$ for each $m=1,\dotsc,M$},
\end{align*}
where $y \mapsto z$ defined by $z = y + h \tgg(y,z,h)$ is a symmetric one-step method applied to \eqref{p1}, $A'' := [ i_1(x) \dotsi i_M(x)] \in \bbR^{d \times M}$ and $A' := [i_1(x') \dotsi i_M(x')] \in \bbR^{d \times M}$.  If we let $\lambda = 2 \mu$ and $A = \tsfrac{1}{2} (A'' + A')$, and eliminate $y$ and $z$, then we can write the method in one line: Given $x\in \bbR^d$ and $h \in [0,\infty)$, compute $x' \in \bbR^d$ and $\lambda \in \bbR^M$ such that
$$
	x' = x + h \tgg\left(x + \tsfrac{1}{2} A'' \lambda, x' - \tsfrac{1}{2} A' \lambda,h\right) + A \lambda \qquad \mbox{and $I_m(x') = I_m(x)$ for each $m$}.
$$
Let $\Psi$ be the implicitly defined map so that $\Psi(x,x',h) = \lambda$ where $\lambda$ satisfies
$$
	\lambda = -h (B^T A)^{-1} B^T \tgg\left(x + \tsfrac{1}{2} A'' \lambda, x' - \tsfrac{1}{2} A' \lambda, h\right)
$$
where $B = [\bii_1 \dotsi \bii_M] \in \bbR^{d \times M}$ for some arbitrarily chosen discrete gradients $\bii_m = \bii_m(x,x')$ of $I_m$.  It is now clear that this method may be written in the form \eqref{eqn0m} if for each $m$ we define $\tii_m : \bbR^d \times \bbR^d \times [0,\infty) \rightarrow \bbR^d$ and $\tff : \bbR^d \times \bbR^d \times [0,\infty) \rightarrow \bbR^d$ such that
$$
	\tii_m(x,x',h) := \tsfrac{1}{2} (i_m(x) + i_m(x')) \qquad \mbox{and} \qquad
	\tff(x,x',h) := \tgg(y,z,h)
$$
for all $x,x' \in \bbR^d$ and $h \in [0,\infty)$ where $y := x + \tsfrac{1}{2} A'' \Psi(x,x',h)$ and $z := x' - \tsfrac{1}{2} A' \Psi(x,x',h)$.

{\noindent\bf Example 3 revisited: Methods of Dahlby, Owren and Yaguchi}.

The multiple integral preserving form of these methods is obtained by replacing the projection matrix $P$ defined earlier in \eqref{eqn8a} with the orthogonal projection matrix
$$
	P = I - B (B^T B)^{-1} B^T
$$
where $B = [\bii_1 \dotsi \bii_M]$ and each $\bii_m = \bii_m(x,x')$ is an arbitrarily chosen discrete gradient of $I_m$.  With this new $P$ the first of the Dahlby et al. methods is (cf. \eqref{eqn8})
$$
	y = \Phi_h(x), \qquad x' = x + P(y-x),
$$
where $\Phi_h$ defines an arbitrary one-step method for solving \eqref{p1}.  This method has the form of a general linear projection method \eqref{eqn0m} if for each $m$ we define $\tii_m : \bbR^d \times \bbR^d \times [0,\infty) \rightarrow \bbR^d$ and $\tff : \bbR^d \times \bbR^d \times [0,\infty) \rightarrow \bbR^d$ such that
$$
	\tii_m(x,x',h) := \bii_m(x,x') \qquad \mbox{and} \qquad \tff(x,x',h) := \frac{\Phi_h(x) - x}{h}
$$
for all $x,x' \in \bbR^d$ and $h \in [0,\infty)$ where each $\bii_m$ is an arbitrary discrete gradient of $I_m$.  The second Dahlby et al. method uses the same choice of $\tii_m$, but a different $\tff$ (defined earlier in \eqref{eqn58}).

\section{Numerical Examples}
\label{sec numerics}


In this section we use a numerical example to provide evidence that the same results for preserving a single first integral, also hold for methods that preserve multiple first integrals.  In particular we will show:
\begin{enumerate}
\item many possible projection directions can be used to define a projection method that preserves the order of accuracy of the underlying method; and
\item the choice of discrete gradient for discrete gradient methods that are equivalent to projection methods does not change the approximate solution in exact arithmetic, but there may be differences in finite precision arithmetic.  
\end{enumerate}

The example we use is Kepler's two-body problem in cartesian coordinates (see e.g. \cite[\S I.2]{HLW} and \cite{dahlby11}).  We will consider the case where either two or three integrals are preserved (the fourth integral is not functionally independent).  

Kepler's two-body problem in the form of \eqref{p1} is
$$
	\frac{\dd}{\dd t} \left[ \fourbyone{x_1}{x_2}{x_3}{x_4} \right] 
	= \left[ \fourbyone{x_3}{x_4}{-x_1 / r^3}{-x_2 / r^3} \right], \qquad t > 0,
$$
where $r = (x_1^2 + x_2^2)^{1/2}$, so that $f(x) := (x_3,x_4,-x_1/r^3, -x_2/r^3)^T$.  This system models two bodies that attract each other, with one body at the origin and the second body at position $(x_1,x_2)$ with velocity (or momentum if the body has mass $1$) $(x_3,x_4)$.  The variable $r$ is the distance between the two bodies.  The exact solution to Kepler's two-body problem preserves four first integrals,
\begin{align*}
	I_1 &:= \tsfrac{1}{2}(x_3^2 + x_4^2) - \tsfrac{1}{r} \\
	I_2 &:= x_1 x_4 - x_3 x_3 \\
	I_3 &:= x_2 x_3^2 - x_1 x_3 x_4 - \tsfrac{x_2}{r} \\
	I_4 &:= x_1 x_4^2 - x_2 x_3 x_4 - \tsfrac{x_1}{r}.
\end{align*}
These integrals are energy, angular momentum, and the two components of the Runge-Lenz-Pauli vector respectively.  As in \cite[p. 12]{HLW}, we use the initial condition 
\begin{equation}
\label{eqn70}
	x_0 = \left( 1-e, 0, 0, \sqrt{\tsfrac{1+e}{1-e}} \right)^T
\end{equation}
for some $e \in [0,1)$ so that the exact solution has period $2 \pi$.  The exact solution can be found by integrating equation (2.10) in \cite[p. 11]{HLW} but we will use a very accurate solution computed with Matlab's ODE45 and very small tolerances as a reference solution in our examples.  

\begin{figure}
\begin{center}
\resizebox{0.48\textwidth}{!}{\input{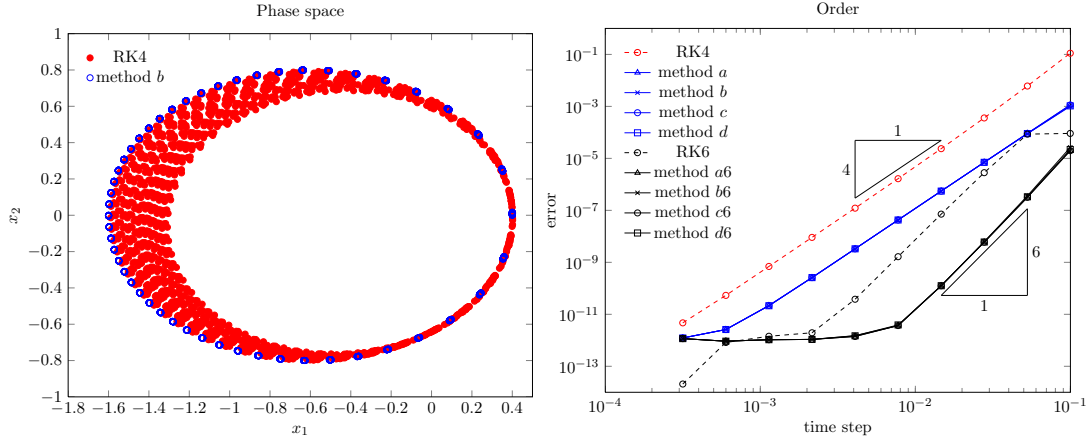}}
\resizebox{0.50\textwidth}{!}{
%
%
%
%
\begin{tikzpicture}

\begin{loglogaxis}[%
view={0}{90},
width=0.7*6.01828521434821in,
height=0.7*4.74667979002625in,
scale only axis,
xmin=0.0001, xmax=0.1,
xminorticks=true,
xlabel={time step},	
ymin=1e-14, ymax=1,
yminorticks=true,
ylabel={error},
title={Order},
legend style={at={(0.03,0.97)},anchor=north west,fill=none,draw=none,align=left}]
\addplot [
color=red,
dashed,
mark=o,
mark options={solid}
]
coordinates{
 (0.000316230575629352,4.68917258365063e-12)(0.000599483380133536,5.35499276389371e-11)(0.00113640537297515,6.95056028817707e-10)(0.00215472747159794,9.11489702937106e-09)(0.00408529603847827,1.21421770884428e-07)(0.00773791293987634,1.65214502620234e-06)(0.014680339502756,2.36304272861647e-05)(0.0278017048990247,0.000359112542326323)(0.0527998765309209,0.0060738846559855)(0.0997331001139617,0.111056186672448) 
};
\addlegendentry{RK4};

\addplot [
color=blue,
solid,
mark=triangle,
mark options={solid}
]
coordinates{
 (0.000316230575629352,1.22931512962328e-12)(0.000599483380133536,2.59482647963335e-12)(0.00113640537297515,2.17574061304477e-11)(0.00215472747159794,2.58533411834404e-10)(0.00408529603847827,3.31099203295807e-09)(0.00773791293987634,4.25944959412968e-08)(0.014680339502756,5.51520423810821e-07)(0.0278017048990247,7.08050790829864e-06)(0.0527998765309209,9.14470884880621e-05)(0.0997331001139617,0.00113091313077064) 
};
\addlegendentry{method $a$};

\addplot [
color=blue,
solid,
mark=x,
mark options={solid}
]
coordinates{
 (0.000316230575629352,1.22931512962328e-12)(0.000599483380133536,2.60184346703845e-12)(0.00113640537297515,2.14097613937203e-11)(0.00215472747159794,2.58176807287344e-10)(0.00408529603847827,3.31044675965787e-09)(0.00773791293987634,4.25700152684867e-08)(0.014680339502756,5.50376501456779e-07)(0.0278017048990247,7.02733698290239e-06)(0.0527998765309209,8.88843653561434e-05)(0.0997331001139617,0.00100330497047721) 
};
\addlegendentry{method $b$};

\addplot [
color=blue,
solid,
mark=o,
mark options={solid}
]
coordinates{
 (0.000316230575629352,1.22931512962328e-12)(0.000599483380133536,2.59482647963335e-12)(0.00113640537297515,2.17574061304477e-11)(0.00215472747159794,2.58533411834404e-10)(0.00408529603847827,3.31099203295807e-09)(0.00773791293987634,4.25944959412968e-08)(0.014680339502756,5.51520424171924e-07)(0.0278017048990247,7.08050779235158e-06)(0.0527998765309209,9.14470216747296e-05)(0.0997331001139617,0.0011308731113714) 
};
\addlegendentry{method $c$};

\addplot [
color=blue,
solid,
mark=square,
mark options={solid}
]
coordinates{
 (0.000316230575629352,1.22931512962328e-12)(0.000599483380133536,2.60050459025463e-12)(0.00113640537297515,2.13753266078192e-11)(0.00215472747159794,2.58237072760058e-10)(0.00408529603847827,3.31059086277801e-09)(0.00773791293987634,4.25780851058463e-08)(0.014680339502756,5.50754253387492e-07)(0.0278017048990247,7.04494741173807e-06)(0.0527998765309209,8.97424177071718e-05)(0.0997331001139617,0.00104776333805077) 
};
\addlegendentry{method $d$};

\addplot [
color=black,
dashed,
mark=o,
mark options={solid}
]
coordinates{
 (0.000316230575629352,2.07771627167895e-14)(0.000599483380133536,8.25344575815679e-13)(0.00113640537297515,1.43068263467194e-12)(0.00215472747159794,1.94902524744376e-12)(0.00408529603847827,3.80316996070351e-11)(0.00773791293987634,1.64783461583419e-09)(0.014680339502756,7.15515691161362e-08)(0.0278017048990247,2.81561856209037e-06)(0.0527998765309209,8.57362806823007e-05)(0.0997331001139617,9.12041373613532e-05) 
};
\addlegendentry{RK6};

\addplot [
color=black,
solid,
mark=triangle,
mark options={solid}
]
coordinates{
 (0.000316230575629352,1.16051063733973e-12)(0.000599483380133536,9.2494788626438e-13)(0.00113640537297515,1.05504392551397e-12)(0.00215472747159794,1.08141355651413e-12)(0.00408529603847827,1.40456979471027e-12)(0.00773791293987634,3.6941988024125e-12)(0.014680339502756,1.25376852650136e-10)(0.0278017048990247,5.96305252433343e-09)(0.0527998765309209,3.17163498391377e-07)(0.0997331001139617,1.97754820079402e-05) 
};
\addlegendentry{method $a6$};

\addplot [
color=black,
solid,
mark=x,
mark options={solid}
]
coordinates{
 (0.000316230575629352,1.16053601208221e-12)(0.000599483380133536,9.24840244713994e-13)(0.00113640537297515,1.05504392551397e-12)(0.00215472747159794,1.08137582296026e-12)(0.00408529603847827,1.51164563235062e-12)(0.00773791293987634,3.6919730996148e-12)(0.014680339502756,1.27076951017163e-10)(0.0278017048990247,6.23103011112367e-09)(0.0527998765309209,3.61318998961e-07)(0.0997331001139617,2.63242157585286e-05) 
};
\addlegendentry{method $b6$};

\addplot [
color=black,
solid,
mark=o,
mark options={solid}
]
coordinates{
 (0.000316230575629352,1.16051063733973e-12)(0.000599483380133536,9.2494788626438e-13)(0.00113640537297515,1.05504392551397e-12)(0.00215472747159794,1.08141355651413e-12)(0.00408529603847827,1.40456979471027e-12)(0.00773791293987634,3.6941988024125e-12)(0.014680339502756,1.25376852650136e-10)(0.0278017048990247,5.96305252433343e-09)(0.0527998765309209,3.17163535892579e-07)(0.0997331001139617,1.97757033631428e-05) 
};
\addlegendentry{method $c6$};

\addplot [
color=black,
solid,
mark=square,
mark options={solid}
]
coordinates{
 (0.000316230575629352,1.16053614486619e-12)(0.000599483380133536,9.24840244713994e-13)(0.00113640537297515,1.05504392551397e-12)(0.00215472747159794,1.08137582296026e-12)(0.00408529603847827,1.48980154706194e-12)(0.00773791293987634,3.89743554941169e-12)(0.014680339502756,1.25887180949446e-10)(0.0278017048990247,6.05929606293965e-09)(0.0527998765309209,3.32742358290279e-07)(0.0997331001139617,2.19374625436206e-05) 
};
\addlegendentry{method $d6$};

\addplot [
color=black,
solid,
forget plot
]
coordinates{
 (0.00408529603847827,4.92958659190598e-05)(0.014680339502756,4.92958659190598e-05) 
};
\addplot [
color=black,
solid,
forget plot
]
coordinates{
 (0.00408529603847827,2.95637985239684e-07)(0.014680339502756,4.92958659190598e-05) 
};
\addplot [
color=black,
solid,
forget plot
]
coordinates{
 (0.00408529603847827,4.92958659190598e-05)(0.00408529603847827,2.95637985239684e-07) 
};
\node[above, inner sep=1mm, text=black]
at (axis cs:0.00774425805446367, 4.92958659190598e-05) {1};
\node[left, inner sep=1mm, text=black]
at (axis cs:0.00408529603847827, 3.8175555635716e-06) { 4};
\addplot [
color=black,
solid,
forget plot
]
coordinates{
 (0.014680339502756,5.32136027314097e-11)(0.0527998765309209,5.32136027314097e-11) 
};
\addplot [
color=black,
solid,
forget plot
]
coordinates{
 (0.014680339502756,5.32136027314097e-11)(0.0527998765309209,1.151870109741e-07) 
};
\addplot [
color=black,
solid,
forget plot
]
coordinates{
 (0.0527998765309209,5.32136027314097e-11)(0.0527998765309209,1.151870109741e-07) 
};
\node[below, inner sep=1mm, text=black]
at (axis cs:0.0278409790269222, 5.32136027314097e-11) {1};
\node[right, inner sep=1mm, text=black]
at (axis cs:0.0527998765309209, 2.47578590386856e-09) { 6};
\end{loglogaxis}
\end{tikzpicture}%
}
\end{center}
\caption{Phase space and order plots for the RK4 and RK6 methods and several projection methods that preserve three integrals.  Methods $a$ -- $d$ and $a6$ -- $d6$ use different projection directions.}
\label{fig1}
\end{figure}

In Figure \ref{fig1} we have computed the solution to the Kepler two-body problem with initial condition \eqref{eqn70} using $e = 0.6$ for several projection methods that differ according to which underlying method is used to define $\tff$ and which projection direction is used to define $\tii$.  We have used the classical explicit $4^{\rm th}$ order and $6^{\rm th}$ order Runge-Kutta methods  (RK4, see e.g. \cite[p. 30]{HLW}, and RK6, see \cite[p. 194]{butcher}, respectively), with coefficients defined by Butcher tableaux:
\renewcommand{\arraystretch}{1.3} 
$$
	\begin{array}{c|cccc}
		\star & & & & \\
		\star & \tsfrac{1}{2} & & &  \\
		\star & 0 & \tsfrac{1}{2} & & \\
		\star & 0 & 0 & 1 & \\ 
		\hline 
		& \tsfrac{1}{6} & \tsfrac{1}{3} & \tsfrac{1}{3} & \tsfrac{1}{6} 
	\end{array}
	\quad \mbox{and} \quad
	\begin{array}{c|ccccccc}
		\star & & & & & &\\
		\star & \tsfrac{1}{3} & & & & &\\
		\star & 0 & \tsfrac{2}{3} & & & & \\
		\star & \tsfrac{1}{12} & \tsfrac{1}{3} & \matminus \tsfrac{1}{12} & & & \\ 
		\star & \tsfrac{25}{48} & \matminus \tsfrac{55}{24} & \tsfrac{35}{48} & \tsfrac{15}{8} & \\
		\star & \tsfrac{3}{20} & \matminus \tsfrac{11}{24} & \matminus \tsfrac{1}{8} & \tsfrac{1}{2} & \tsfrac{1}{10} \\
		\star & \matminus \tsfrac{261}{260} & \tsfrac{33}{13} & \tsfrac{43}{156} & \matminus \tsfrac{118}{39} & \tsfrac{32}{195} & \tsfrac{80}{39} \\ 
		\hline
		& \tsfrac{13}{200} & 0 & \tsfrac{11}{40} & \tsfrac{11}{40} & \tsfrac{4}{25} & \tsfrac{4}{25} & \tsfrac{13}{200}
	\end{array}.
$$
Define $\tff$ using these coefficients by $\tff := h \sum_{i=1}^{s} b_i k_i$ where $b_i$ and $k_i$ are defined as in \cite[eq. (1.4) on p. 29]{HLW}.  Note that we do not require the coefficients indicated by $\star$ in the tableaux because we are solving an autonomous ODE.  

We then define projection methods $a$, $b$, $c$, $d$, $a6$, $b6$, $c6$ and $d6$ for preserving three integrals ($M = 3$) by \eqref{eqn0m} by choosing RK4 to define $\tff$ for methods $a - d$ and RK6 to define $\tff$ for methods $a6 - d6$.  The projection directions for these methods are defined by
\begin{align*}
	\mbox{methods $a$ and $a6$:} &&& \tii^m := i^m(x'), \\
	\mbox{methods $b$ and $b6$:} &&& \tii^m := i^m(x), \\
	\mbox{methods $c$ and $c6$:} &&& \tii^m := i^m(y) \mbox{ where $y = x + h \tff(x,h)$,} \\
	\mbox{methods $d$ and $d6$:} &&& \tii^m := \tsfrac{1}{2} (i^m(x) + i^m(x')).
\end{align*}
We modify \eqref{eqn0m} in our computations to prevent the value of $I_m$ drifting due to finite precision arithmetic.  Instead of requiring $I_m(x') = I_m(x)$ at each time step, we require that $I_m(x') = I_m(x_0)$ for each $m$.

In the phase space plot of Figure \ref{fig1} (left) we see that method $b$ does indeed keep the approximate solution on the ellipse while the RK4 approximate solution drifts away from the ellipse.  In this plot we computed the solution up to a final time of $t=50 \pi$ ($25$ periods) with $h = \tsfrac{2 \pi}{50}$.

In the order plot of Figure \ref{fig1} (right) we see the more important result that methods $a$ -- $d$ all seem to preserve the $4^{\rm th}$ order convergence of the RK4 method and methods $a6$ -- $d6$ all seem to preserve the $6^{\rm th}$ order convergence of the RK6 method.  In fact, we see that the different choices for the projection direction seem to make very little difference to the error because the errors are essentially the same size for methods $a$ -- $d$ and methods $a6$ -- $d6$ respectively (the lines in the plot overlay each other).  For this plot we computed up to a final time of $t=2 \pi$ (only $1$ period) for a range of step sizes in $[10^{-3.5}, 10^{-1}]$.

\begin{figure}
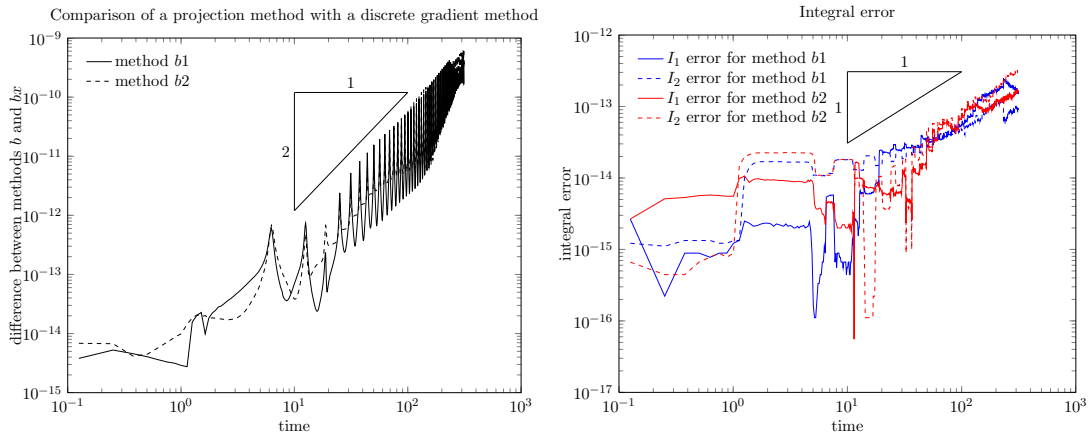

\begin{center}
\resizebox{0.49\textwidth}{!}{\input{fig2-comparison.tikz}}
\resizebox{0.49\textwidth}{!}{\input{fig2-interr.tikz}}
\end{center}
\caption{Left: A plot of the difference between method $b$ to methods $b1$ and $b2$ for increasing time. Right: A plot of the integral error of $I_1$ and $I_2$ for methods $b1$ and $b2$ as time increases.}
\label{fig2}
\end{figure}

In Figure \ref{fig2} we have computed the solution to the Kepler two-body problem using three methods that are equivalent in exact arithmetic according to our theory - a projection method, and two discrete gradient methods defined using different choices of discrete gradient.  To avoid making the discrete gradient methods overly complicated we only consider the case when two integrals, $I_1$ and $I_2$, are preserved.  The projection method we use is method $b$ as defined above except now we only preserve two integrals.  We define methods $b1$ and $b2$ by \eqref{eqn61} and \eqref{eqn62} with $M=2$, the same choices of $\tff$ and $\tii^m$ as for method $b$, and if $\bjj^m$ denotes the coordinate increment discrete gradient of $I_m$(see \cite{abeitoh} or \cite{MQR99}) then $\bii^m$ in methods $b1$ and $b2$ is defined as
\begin{align*}
	\mbox{method $b1$:} &&& \bii^m = \bii^m(x,x') = \bjj^m(x,x'), \\
	\mbox{method $b2$:} &&& \bii^m = \bii^m(x,x') = \tsfrac{1}{2} (\bjj^m(x,x') + \bjj^m(x',x)).
\end{align*}
In exact arithmetic, according to Theorems \ref{thm equiv1m} and \ref{thm equiv3} these methods should be the same.  However, in finite precision arithmetic we notice some small differences.


In the left plot of Figure \ref{fig2} we have  compared methods $b1$ and $b2$ with method $b$ for increasing time.  Since computations are done in finite precision arithmetic and the nonlinear systems at each time step are only solved to a tolerance of $10^{-14}$ we expect to see that the differences between these methods grows linearly with time.  Perhaps surprisingly we actually see quadratic growth in the difference between these methods.  In the right plot of Figure \ref{fig2} we have plotted the error of $I_1$ and $I_2$ for the approximate solution as time increases for methods $b1$ and $b2$ (method $b$ is constructed to keep the integral error below $10^{-14}$, the tolerance that we solve the nonlinear systems at each time step).  We see linear growth of the integral errors as time increases, as expected.  The plots in Figure \ref{fig2} used the same initial condition as above, a time step of $h = \tsfrac{2 \pi}{50}$ and a final time of $t = 100 \pi$ (50 periods).

\section{Conclusions}
\label{sec conclusion}

In this paper we determined the relationship between linear projection methods and discrete gradient methods for ODEs with conserved first integrals.  A consequence of our theory is that each linear projection method is equivalent to a class of discrete gradient methods.  A further consequence when there is only one first integral to preserve is that we can use theory from discrete gradient methods to prove results about projection methods.  We have shown that under only mild conditions on the continuity and consistency of the projection direction we obtain a projection method that has a well-defined approximate solution provided the time step is sufficiently small, and for arbitrary $p \in \mathbb{N}$ also preserves the order of accuracy of an underlying method of order $p$.  Moreover, the condition on the projection direction does not depend on $p$.  For the multiple first integral case we rely on numerical experiments to confirm that similar results appear to also hold in this case.  

\section*{Acknowledgements}

This research was supported by the Australian Research Council and the Marie Curie International Research Staff Exchange Scheme within the 7th European Community Framework Programme.


\bibliographystyle{plain}
\bibliography{paper7bib}

\begin{appendix}
\section{Example 2 continued from \S \ref{sec ex2}}
Continuing on from the end of \S \ref{sec ex2} we seek to show that $\tgg$ satisfying Assumption \ref{assumption1} is sufficient to ensure that both $\tii$ and $\tff$ defined by \eqref{eqn53} satisfy Assumptions \ref{assumption1} and \ref{assumption2} respectively.  We begin by verifying that $\tii$ defined by \eqref{eqn53} satisfies Assumption \ref{assumption2} for some choice of $R$, $L$ and $H$.

\begin{lemma}
\label{lem tii3}
For a bounded set $B \subset \bbR^d$, suppose that $i$ satisfies Assumption \ref{assumption3} for some positive constants $R$ and $L_i$, and let $H$ be an arbitrary positive constant.  Define $L := \tsfrac{1}{2} L_i$.  Then $\tii : \bbR^d \times \bbR^d \times [0,\infty) \rightarrow \bbR^d$ defined by $\tii(x,x',h) := \tsfrac{1}{2} (i(x) + i(x'))$ for all $x,x' \in \bbR^d$ and $h \in [0,\infty)$ satisfies Assumption \ref{assumption2} for $R$, $L$ and $H$.
\end{lemma}

\begin{proof}
Since $i$ is locally Lipschitz continuous, $L_i$ exists for any given $R > 0$.  Let $x \in B$, $u,v,w \in B_R(x)$ and $h \in [0,\infty)$.  Trivially, we get $\tii(x,x,0) =i(x)$ and $|\tii(x,x,h) - \tii(x,x,0)| = 0$.  Using Assumption \ref{assumption3} for $i$ we also get
$$
	|\tii(u,v,h) - \tii(w,v,h)| = \tsfrac{1}{2} |i(u) - i(w)| \leq \tsfrac{1}{2} L_i |u-w| = L |u-w|.
$$
Similarly, we get $|\tii(u,v,h) - \tii(u,w,h)| \leq L |v-w|$. 
\end{proof}

Verifying that $\tff$ defined by \eqref{eqn53} satisfies Assumption \ref{assumption1} requires a much lengthier argument.  We first prove the following lemma to describe the properties of $\Psi$.

\begin{lemma}
\label{lem tff2}
For a bounded set $B \subset \bbR^d$, let $C_1$ be the constant from \eqref{eqn2}, and suppose that
\begin{enumerate}
\item $\tgg : \bbR^d \times \bbR^d \times [0,\infty) \rightarrow \bbR^d$ satisfies Assumption \ref{assumption1} for positive constants $R_g$, $L_g$ and $H_g$,
\item $\tii : \bbR^d \times \bbR^d \times [0,\infty) \rightarrow \bbR^d$ satisfies Assumption \ref{assumption2} for positive constants $R_{\tii}$, $L_{\tii}$ and $H_{\tii}$,
\item $\bii$ is a discrete gradient of $I$ satisfying Assumption \ref{assumption2} for $R_{\tii}$ and $L_{\tii}$, and
\item $i$ satisfies Assumption \ref{assumption3} for $R_g$ and a positive constant $L_i$.
\end{enumerate}
Define 
\begin{align*}
	R_\lambda := \max \{ 2R_g, R_{\tii}, 12 L_{\tii}, L_i, 4 L_g \},
	\quad \mbox{and} \quad 
	H_\lambda := \min \left\{ H_g, H_{\tii}, \tsfrac{1}{9 L_{\tii}}, \tsfrac{3}{28 L_g}, \tsfrac{3}{7(C_1+2) R_\lambda} \right\}.
\end{align*}
For any $x \in B$ such that $i(x) \neq 0$, any $u,v \in B_{R_\lambda}(x)$ and any $h \in [0,H_\lambda)$, there exists a unique $\lambda = \Psi(u,v,h) \in \bbR$ such that $|\lambda| \leq \tsfrac{1}{R_\lambda}$, satisfying
$$
	\lambda = - h \frac{ \tgg\left(u + \tsfrac{\lambda}{2} i(u), v - \tsfrac{\lambda}{2} i(v),h\right) \cdot \bii(u,v)}{\tii(u,v,h) \cdot \bii(u,v)}.
$$
Moreover, 
\begin{equation}
\label{eqn9s0}
	|\Psi(u,v,h)| \leq \tsfrac{7}{3} (C_1+2) h,
\end{equation}
and if $w \in B_{R_\lambda}(x)$ then 
\begin{equation}
\label{eqn9s1}
\begin{split}
	|\Psi(u,v,h) - \Psi(w,v,h)| &\leq \tsfrac{5}{4} \tsfrac{1}{|i(x)|} |u-w|,  \\
	|\Psi(u,v,h) - \Psi(u,w,h)| &\leq \tsfrac{5}{4} \tsfrac{1}{|i(x)|} |v-w|.  
\end{split}
\end{equation}
If $x \in B$ such that $i(x) = 0$, then define $\lambda = \Psi(u,v,h) :=0$ for all $u,v \in B_{R_\lambda}(x)$ and $h \in [0,H_\lambda)$.
\end{lemma}

\begin{proof}
Note that since $i$ is locally Lipschitz continuous, $L_i$ exists for any $R_g > 0$.  Fix $x \in B$ such that $i(x) \neq 0$, $u,v \in B_{R_\lambda}(x)$ and $h \in [0,H_\lambda)$.  To get the result we will use Theorem \ref{thm banach}.  Define $X := \lbrace \gamma \in \bbR : |\gamma| \leq \tsfrac{1}{R_\lambda} \rbrace$ (which with the Euclidean norm $|\cdot|$ is a non-empty complete Metric space) and $T : X \rightarrow \bbR$ by 
$$
	T(\gamma):= - h \frac{ \tgg\left(u + \tsfrac{\gamma}{2} i(u), v - \tsfrac{\gamma}{2} i(v),h\right) \cdot \bii(u,v)}{\tii(u,v,h) \cdot \bii(u,v)} \qquad \mbox{for each $\gamma \in X$}.
$$
To apply Theorem \ref{thm banach} we must show that $T(\gamma) \in X$ for each $\gamma \in X$ and that $T$ is a contraction on $X$.  

Using Assumption \ref{assumption3} for $i$ and $R_\lambda \geq L_i$ we get the following useful inequality,
\begin{equation}
\label{eqn11s}
	|i(u)| \leq |i(x)| + L_i |u-x| \leq \left( 1 + \tsfrac{L_i}{R_\lambda} \right) |i(x)| \leq 2 |i(x)|.
\end{equation}
Fix $\gamma \in X$ and define $y := u + \tsfrac{\gamma}{2} i(u)$ and $z := v - \tsfrac{\gamma}{2} i(v)$.  Using \eqref{eqn11s}, $|\gamma| \leq \tsfrac{1}{R_\lambda}$ and $u \in B_{R_\lambda}(x)$ we get
$$
	|y - x| \leq |u-x| + \tsfrac{|\gamma|}{2} |i(u)| 
	\leq |u-x| + |\gamma| |i(x)| 
	\leq \left( \tsfrac{1}{R_\lambda} + \tsfrac{1}{R_\lambda} \right) |i(x)| 
	= \tsfrac{|i(x)|}{R_\lambda/2}.
$$
Hence $y \in B_{R_\lambda/2}(x)$, and since $R_\lambda/2 \geq R_g$ we have $y \in B_{R_\lambda/2} \subset B_{R_g}(x)$.  Similarly, $z \in B_{R_\lambda/2} \subset B_{R_g}(x)$.  Using this, Assumption \ref{assumption1} for $\tgg$, \eqref{eqn2}, $R_\lambda \geq 4 L_g$ and $h \leq \tsfrac{3}{28L_g} < \tsfrac{1}{L_g}$ we get
\begin{equation}
\label{eqn8s}
\begin{split}
	|\tgg(y,z,h)| 
	&\leq |\tgg(x,x,0)| + L_g (|y-x| + |z-x| + h |i(x)|) \\
	&\leq \left(C_1 + \tsfrac{2L_g}{R_\lambda} + \tsfrac{2L_g}{R_\lambda} + L_g h \right) |i(x)| 
	\leq (C_1 + 2) |i(x)|.
\end{split}
\end{equation}
Using Assumption \ref{assumption2} for $\bii$, and $R_\lambda \geq 12 L_{\tii}$ we get
\begin{equation}
\label{eqn7s}
	|\bii(u,v)| \leq |\bii(x,x)| + L_{\tii} (|u-x| + |v-x|) \leq \left( 1 + \tsfrac{2L_{\tii}}{R_\lambda} \right) |i(x)| \leq \tsfrac{7}{6} |i(x)|.
\end{equation}
Using Assumpiton \ref{assumption2} for $\tii$ and $\bii$, \eqref{eqn7s}, $R_\lambda \geq 12 L_{\tii}$ and $h \leq \tsfrac{1}{9 L_{\tii}}$ we get
\begin{equation}
\label{eqn9s}
\begin{split}
	\tii(u,v,h) \cdot \bii(u,v) 
	&= |i(x)|^2 + \Bigl( [\tii(u,v,h) - \tii(x,v,h)] + [\tii(x,v,h) - \tii(x,x,h)] \\
	&\qquad  + [\tii(x,x,h) - \tii(x,x,0)] \Bigr) \cdot \bii(u,v) \\
	&\qquad + i(x) \cdot \Bigl( [\bii(u,v) - \bii(x,v)] + [\bii(x,v) - \bii(x,x)] \Bigr) \\
	&\geq |i(x)|^2 - L_{\tii} (|u-x| + |v-x| + h |i(x)|) \tsfrac{7}{6} |i(x)| \\
	&\qquad - |i(x)| L_{\tii} (|u-x| + |v-x|) \\
	&\geq |i(x)|^2 \left( 1 - L_{\tii} (\tsfrac{2}{R_\lambda} + h)\tsfrac{7}{6}  - L_{\tii} \tsfrac{2}{R_\lambda} \right) \\
	&\geq |i(x)|^2 \left( 1 - (\tsfrac{1}{6} + \tsfrac{1}{9}) \tsfrac{7}{6} - \tsfrac{1}{6} \right)
	= \tsfrac{55}{108} |i(x)|^2 > \tsfrac{1}{2} |i(x)|^2.
\end{split}
\end{equation}
We can now show that $|T(\gamma)| \leq \tsfrac{1}{R_\lambda}$.  Using \eqref{eqn8s}, \eqref{eqn7s}, \eqref{eqn9s} and $h \leq \tsfrac{3}{7(C_1 + 2) R_\lambda}$ we get
\begin{equation}
\label{eqn10s}
	|T(\gamma)| 
	= h \tsfrac{|\tgg(y,z,h) \cdot \bii(u,v)|}{|\tii(u,v,h) \cdot \bii(u,v) |} 
	\leq \tsfrac{7}{3} (C_1 + 2) h 
	\leq \tsfrac{1}{R_\lambda}.
\end{equation}
Hence $T(\gamma) \in X$ and so $T: X \rightarrow X$.  It remains to show that $T$ is a contraction.  Let $\delta \in X$ and define $y' := u + \tsfrac{\delta}{2} i(u)$ and $z' := v - \tsfrac{\delta}{2} i(v)$.  As above we have $y',z' \in B_{R_\lambda/2}(x) \subset B_{R_g}(x)$.  Using \eqref{eqn9s}, \eqref{eqn7s}, Assumption \ref{assumption1} for $\tgg$, \eqref{eqn11s} (which also holds for $|i(v)|$) and $h \leq \tsfrac{3}{28L_g}$ we get
\begin{align*}
	|T(\gamma) - T(\delta)| 
	&= h \left| \tsfrac{ [\tgg(y,z,h) - \tgg(y',z',h)] \cdot \bii(u,v) }{  \tii(u,v,h) \cdot \bii(u,v) } \right| \\
	&\leq \tsfrac{7 h}{3 |i(x)|} |\tgg(y,z,h) - \tgg(y',z',h)| \\
	&\leq \tsfrac{7 h L_g}{3 |i(x)|} ( |y-y'| + |z-z'| ) \\
	&\leq \tsfrac{7 h L_g}{6 |i(x)|} |\gamma - \delta| ( |i(u)| + |i(v)| ) \\
	&\leq \tsfrac{7 h L_g}{3} |\gamma - \delta| 
	\leq \tsfrac{1}{4} |\gamma - \delta|,
\end{align*}
so $T$ is a contraction.  Therefore, applying Theorem \ref{thm banach}, there exists a unique $\lambda \in X$ such that $T(\lambda) = \lambda$.  Define $\Psi(u,v,h):= \lambda$.  Inequality \eqref{eqn9s0} then follows from \eqref{eqn10s}.

Now let $w \in B_{R_\lambda}(x)$ and define $\lambda := \Psi(u,v,h)$, $\gamma := \Psi(w,v,h)$ and
\begin{align*}
	y &:= u + \tsfrac{\lambda}{2} i(u), & y' &:= w + \tsfrac{\gamma}{2} i(w), \\
	z &:= v - \tsfrac{\lambda}{2} i(v), & z' &:= v - \tsfrac{\gamma}{2} i(v).
\end{align*}
Using Assumption \ref{assumption3} for $i$, $|\lambda| \leq \tsfrac{1}{R_\lambda}$, \eqref{eqn11s} and $R_\lambda \geq L_i$ we get
\begin{align*}
	|y-y'| &\leq |u-w| + \tsfrac{|\lambda|}{2} |i(u) - i(w)| + \tsfrac{|i(w)|}{2} |\lambda - \gamma| \\
	&\leq \left( 1 + \tsfrac{L_i}{2 R_\lambda} \right) |u-w| + |i(x)| |\lambda - \gamma| \\
	&\leq \tsfrac{3}{2} |u-w| + |i(x)| |\lambda - \gamma|.
\end{align*}
Likewise, $|z - z'| \leq \tsfrac{|i(v)|}{2} |\lambda - \gamma| \leq |i(x)| |\lambda - \gamma|$.  Using these two inequalities, together with Assumption \ref{assumption1} for $\tgg$, and noting that $y,y',z,z' \in B_{R_\lambda/2}(x) \subset B_{R_g}(x)$ we get
\begin{equation}
\label{eqn12s}
	|\tgg(y,z,h) - \tgg(y',z',h)| \leq L_g ( |y-y'| + |z-z'| ) \leq L_g \left( \tsfrac{3}{2} |u-w| + 2 |i(x)| |\lambda - \gamma| \right).
\end{equation}
Using Assumption \ref{assumption2} for $\tii$, $u,v \in B_{R_\lambda}(x)$, $R_\lambda \geq 12 L_{\tii}$ and $h \leq \tsfrac{1}{9 L_{\tii}}$ we get
\begin{equation}
\label{eqn13s}
\begin{split}
	|\tii(u,v,h)| 
	&\leq |i(x)| + L_{\tii} ( |u-x| + |v-x| + h |i(x)| ) \\
	&\leq \left(1 + \tsfrac{2 L_{\tii}}{R_\lambda} + L_{\tii} h \right) |i(x)|\\ 
	&\leq (1 + \tsfrac{1}{6} + \tsfrac{1}{9})|i(x)| = \tsfrac{23}{18} |i(x)| < \tsfrac{4}{3} |i(x)|.
\end{split}
\end{equation}
The same inequality holds for $|\tii(w,v,h)|$.  Using \eqref{eqn7s} and \eqref{eqn9s} (which also hold with $u$ replaced by $w$), \eqref{eqn13s}, and Assumption \ref{assumption2} for $\tii$ and $\bii$ we get
\begin{equation}
\label{eqn14s}
\begin{split}
	\Bigl| \tsfrac{1}{\tii(u,v,h)\cdot \bii(u,v)} - &\tsfrac{1}{\tii(w,v,h)\cdot \bii(w,v)} \Bigr|
	= \left| \tsfrac{\tii(w,v,h)\cdot \bii(w,v) - \tii(u,v,h)\cdot \bii(u,v)}{(\tii(u,v,h)\cdot \bii(u,v)) (\tii(w,v,h)\cdot \bii(w,v))} \right| \\
	&= \left| \tsfrac{[\tii(w,v,h) - \tii(u,v,h)] \cdot \bii(w,v) + \tii(u,v,h)\cdot [\bii(w,v) - \bii(u,v)]}{(\tii(u,v,h)\cdot \bii(u,v)) (\tii(w,v,h)\cdot \bii(w,v))} \right| \\
	&\leq \tsfrac{4}{|i(x)|^4} \left( |\tii(u,v,h) - \tii(w,v,h)| \tsfrac{7}{6} |i(x)| + \tsfrac{4}{3} |i(x)| |\bii(u,v) - \bii(w,v)|  \right) \\
	&\leq \tsfrac{4}{|i(x)|^3} \left( \tsfrac{7}{6} L_{\tii} |u-w| + \tsfrac{4}{3} L_{\tii} |u-w| \right) 
	= \tsfrac{10 L_{\tii}}{|i(x)|^3} |u-w|.
\end{split}
\end{equation}
Now using \eqref{eqn9s}, \eqref{eqn7s}, \eqref{eqn8s}, \eqref{eqn12s}, Assumption \ref{assumption2} for $\bii$, \eqref{eqn14s}, $h \leq \min \lbrace \tsfrac{3}{28L_g}, \tsfrac{1}{7(C_1+2) R_\lambda} \rbrace$ and $R_\lambda \geq 12 L_{\tii}$ we get
\begin{equation}
\label{eqn15s}
\begin{split}
	|\lambda - \gamma| 
	=& h \left| \tsfrac{\tgg(y,z,h)\cdot\bii(u,v)}{\tii(u,v,h)\cdot\bii(u,v)} - \tsfrac{\tgg(y',z',h)\cdot\bii(w,v)}{\tii(w,v,h)\cdot\bii(w,v)}  \right|\\
	=& h \biggl| 
		 {\textstyle \frac{[\tgg(y,z,h) - \tgg(y',z',h] \cdot \bii(u,v)}{\tii(u,v,h)\cdot\bii(u,v)}  
	 + \frac{\tgg(y',z',h) \cdot [ \bii(u,v) - \bii(w,v)]}{\tii(u,v,h)\cdot\bii(u,v)}  }\\
	& \qquad + {\scriptstyle \tgg(y',z',h) \cdot \bii(w,v)} \left( \tsfrac{1}{\tii(u,v,h)\cdot \bii(u,v)} - \tsfrac{1}{\tii(w,v,h)\cdot \bii(w,v)} \right)
		 \biggr| \\
	\leq& h \biggl( {\scriptstyle \tsfrac{2}{|i(x)|^2} \tsfrac{7}{6}  |i(x)| | \tgg(y,z,h) - \tgg(y',z',h)| 
		+ \tsfrac{2}{|i(x)|^2} (C_1 + 2) |i(x)| |\bii(u,v) - \bii(w,v)|} \\
	& \qquad  + {\scriptstyle (C_1 + 2) |i(x)| \tsfrac{7}{6} |i(x)| \left| \tsfrac{1}{\tii(u,v,h)\cdot \bii(u,v)} - \tsfrac{1}{\tii(w,v,h)\cdot \bii(w,v)} \right| }
	\biggr) \\
	=& h \biggl( {\scriptstyle \tsfrac{7}{3 |i(x)|}  |\tgg(y,z,h) - \tgg(y',z',h)| 
	+ \tsfrac{2(C_1 + 2)}{|i(x)|}   |\bii(u,v) - \bii(w,v)|} \\
	& \qquad + \tsfrac{7(C_1 + 2)|i(x)|^2}{6}  \left| \tsfrac{1}{\tii(u,v,h)\cdot \bii(u,v)} - \tsfrac{1}{\tii(w,v,h)\cdot \bii(w,v)} \right| 
	\biggr) \\
	\leq& h \biggl( {\scriptstyle \tsfrac{7L_g}{3 |i(x)|}  \left( \tsfrac{3}{2} |u-w| + 2 |i(x)| |\lambda - \gamma| \right) 
	+ \tsfrac{2(C_1 + 2)}{|i(x)|} L_{\tii} |u-w| 
	+  \tsfrac{7(C_1 + 2)|i(x)|^2}{6}  \tsfrac{10L_{\tii}}{|i(x)|^3} |u-w| \biggr) }\\
	=& \left( \tsfrac{7L_g h}{2} + \tsfrac{41L_{\tii} (C_1+2) h}{3} \right) \tsfrac{1}{|i(x)|}|u-w| + \tsfrac{14 L_g }{3} h |\lambda - \gamma| \\
	\leq& \left( \tsfrac{3}{8} + \tsfrac{41L_{\tii}}{21 R_\lambda} \right) \tsfrac{1}{|i(x)|}|u-w| + \tsfrac{1}{2} |\lambda - \gamma| \\
	\leq& \left( \tsfrac{3}{8} + \tsfrac{1}{4} \right) \tsfrac{1}{|i(x)|}|u-w| + \tsfrac{1}{2} |\lambda - \gamma| \\
	=& \tsfrac{5}{8} \tsfrac{1}{|i(x)|}|u-w| + \tsfrac{1}{2} |\lambda - \gamma|.
\end{split}
\end{equation}
It then follows that $|\Psi(u,v,h) - \Psi(w,v,h)| \leq \tsfrac{5}{4} \tsfrac{1}{|i(x)|}|u-w|$.  The second inequality in \eqref{eqn9s1} is derived using a similar argument.
\end{proof}

Using Lemma \ref{lem tff2} we can easily derive the following lemma to ensure that $\tff$ for the symmetric projection method (defined by \eqref{eqn53}) satisfies Assumption \ref{assumption1}.

\begin{lemma}
\label{lem tff2a}
For a bounded set $B \subset \bbR^d$, define $C_1$, $\tgg$, $R_g$, $L_g$, $H_g$, $\tii$, $R_{\tii}$, $L_{\tii}$, $H_{\tii}$, $\bii$, $L_{i}$, $R_\lambda$ and $H_\lambda$ as in Lemma \ref{lem tff2}.  In addition, define
$$
	L := \max \left\lbrace 4 L_g , \tsfrac{(7 C_1+17)L_g}{3} \right\rbrace.
$$
Then $\tff : \bbR^d \times \bbR^d \times [0,\infty) \rightarrow \bbR^d$, defined as in \eqref{eqn53},
$$
	\tff(x,x',h) := \tgg \left(x + \tsfrac{\Psi(x,x',h)}{2} i(x), x' - \tsfrac{\Psi(x,x',h)}{2} i(x'),h \right),
$$
for all $x,x' \in \bbR^d$ and $h \in [0,\infty)$ satisfies Assumption \ref{assumption1} for $R_\lambda$, $L$ and $H_\lambda$.
\end{lemma}

\begin{proof}
Since the assumptions of Lemma \ref{lem tff2a} are the same as Lemma \ref{lem tff2} all of the results in Lemma \ref{lem tff2} hold here.  Fix $x \in B$, $u,v,w \in B_{R_\lambda}(x)$ and $h \in [0,H_\lambda)$.  Since $\Psi(x,x,0) = 0$ (by Lemma \ref{lem tff2}) it follows from Assumption \ref{assumption1} for $\tgg$ that $\tff(x,x,0) = \tgg(x,x,0) = f(x)$.  Using Lemma \ref{lem tff2} we can define $\lambda := \Psi(u,v,h)$ and $\gamma := \Psi(w,v,h)$, and
\begin{align*}
	y &:= u + \tsfrac{\lambda}{2} i(u), & y' &:= w + \tsfrac{\gamma}{2} i(w), \\
	z &:= v - \tsfrac{\lambda}{2} i(v), & z' &:= v - \tsfrac{\gamma}{2} i(v).
\end{align*}
From \eqref{eqn12s} in the proof of Lemma \ref{lem tff2}, and \eqref{eqn9s1} we get 
\begin{align*}
	|\tff(u,v,h) - \tff(w,v,h)| 
	&= |\tgg(y,z,h) - \tgg(y',z',h)| \\
	&\leq L_g \left( \tsfrac{3}{2} |u-w| + 2 |i(x)| |\lambda -\gamma| \right) \\
	&\leq L_g \left( \tsfrac{3}{2} + \tsfrac{5}{2} \right) |u-w| 
	= 4 L_g |u-w| \leq L |u-w|.
\end{align*}
Similarly, we can show that $|\tff(u,v,h) - \tff(u,w,h)| \leq L|v-w|$.  

Now define $\lambda := \Psi(x,x,h)$, $y := x + \tsfrac{\lambda}{2} i(x)$ and $z := x - \tsfrac{\lambda}{2} i(x)$.  Using $\Psi(x,x,0) = 0$, Assumption \ref{assumption1} for $\tgg$ and \eqref{eqn9s0} we get
\begin{align*}
	|\tff(x,x,h) - \tff(x,x,0)|
	&= |\tgg(y,z,h) - \tgg(x,x,0)| \\
	&\leq L_g ( |y-x| + |z-x| + h |i(x)| ) 
	= L_g (|\lambda| + h) |i(x)| \\
	&\leq L_g \left( \tsfrac{7}{3}(C_1 + 2) + 1 \right) h |i(x)| 
	= \tsfrac{(7C_1 + 17)L_g}{3}  h |i(x)| \leq L h |i(x)|.
\end{align*}
Hence $\tff$ satisfies Assumption \ref{assumption1} for $R_\lambda$, $L$ and $H_\lambda$.
\end{proof}
\end{appendix}

\end{document}